\documentclass[11pt,a4paper,final,twoside]{article}
\usepackage[T1]{fontenc}
\usepackage[utf8]{inputenc}
\usepackage[english]{babel}
\usepackage{csquotes}
\usepackage{natbib}
\usepackage[centertags]{amsmath}
\usepackage{amsfonts}
\usepackage{amssymb}
\usepackage{eqparbox}
\usepackage{amsthm}
\usepackage{newlfont}
\usepackage{fancyhdr}
\usepackage{amscd}
\usepackage{scalerel}
\usepackage{graphicx}
\usepackage{subfigure}
\usepackage{afterpage}
\usepackage{relsize}
\usepackage{sectsty}
\usepackage{siunitx}
\usepackage{etoolbox}
\usepackage{caption}
\usepackage{array}
\usepackage{booktabs}
\usepackage[table]{xcolor}
\usepackage{xcolor,colortbl}
\usepackage{colortbl,multirow,hhline}
\definecolor{blue(ncs)}{rgb}{0.0, 0.53, 0.74}
\usepackage[colorlinks=true, urlcolor=blue(ncs), linkcolor=blue(ncs), citecolor=blue(ncs)]{hyperref}
\usepackage{latexsym,epstopdf,float}
\usepackage{euscript,pifont,mathrsfs,latexsym,graphicx}
\usepackage{mathrsfs}
\usepackage{etoolbox}
\usepackage[atend]{bookmark}
\bookmarksetup{
  open,
  openlevel=2,
  numbered,
  addtohook={%
    \ifnum\bookmarkget{level}>1 %
      \DisableBookmarkNumbering
    \fi
  },
}
\makeatletter
\newcommand*{\DisableBookmarkNumbering}{%
  \let\numberline\@gobble
}
\makeatother
\patchcmd\thebibliography{\labelsep}{\labelsep\itemsep=0pt\parsep=0pt\relax}{}{\typeout{Couldn't patch the command}}

\theoremstyle{plain}
\newtheorem{thm}{Theorem}[section]
\newtheorem{prop}[thm]{Proposition}
\newtheorem{lem}[thm]{Lemma}
\newtheorem{cor}[thm]{Corollary}
\theoremstyle{definition}
\newtheorem{defn}{Definition}[section]
\newtheorem*{assum}{Assumptions}
\theoremstyle{remark}
\newtheorem{rem}{Remark}[section]
\newtheorem{exa}[rem]{Example}

 \numberwithin{equation}{section}

 \DeclareMathOperator*{\medcup}{\mathbin{\scalebox{1.1}{\ensuremath{\bigcup}}}}
 
 \DeclareMathOperator*{\medPl}{\mathbin{\scalebox{0.85}{\ensuremath{\Big(}}}}
 \DeclareMathOperator*{\medPr}{\mathbin{\scalebox{0.85}{\ensuremath{\Big)}}}}
 \DeclareMathOperator*{\medBl}{\mathbin{\scalebox{0.85}{\ensuremath{\Big\{}}}}
 \DeclareMathOperator*{\medBr}{\mathbin{\scalebox{0.85}{\ensuremath{\Big\}}}}}

 \newcommand*{\sOmega}{\mathbin{\scalebox{.5}{\ensuremath{\Omega}}}}

\makeatletter
\newcommand{\pushright}[1]{\ifmeasuring@#1\else\omit\hfill$\displaystyle#1$\fi\ignorespaces}
\newcommand{\pushleft}[1]{\ifmeasuring@#1\else\omit$\displaystyle#1$\hfill\fi\ignorespaces}
\makeatother

\begin{document}
\pdfbookmark[0]{Robust optimality and duality for composite uncertain multiobjective optimization in Asplund spaces with its applications}{frontmatter}
\title{Robust optimality and duality for composite uncertain multiobjective optimization in Asplund spaces with its applications}
\author{Maryam Saadati$^*$, Morteza Oveisiha}
\date{}
\maketitle
\begin{center}
Department of Pure Mathematics, Faculty of Science, Imam Khomeini International University, P.O. Box 34149-16818, Qazvin, Iran.
\\
E-mail: m.saadati@edu.ikiu.ac.ir, oveisiha@sci.ikiu.ac.ir
\end{center}
\pdfbookmark[section]{\abstractname}{abstract}
\begin{abstract}
This article is devoted to investigate a nonsmooth/nonconvex uncertain multiobjective optimization problem with composition fields ($(\hyperlink{CUP}{\mathrm{CUP}})$ for brevity) over arbitrary Asplund spaces. Employing some advanced techniques of variational analysis and generalized differentiation, we establish necessary optimality conditions for weakly robust efficient solutions of $(\hyperlink{CUP}{\mathrm{CUP}})$ in terms of the limiting subdifferential. Sufficient conditions for the existence of (weakly) robust efficient solutions to such a problem are also driven under the new concept of pseudo-quasi convexity for composite functions. We formulate a Mond-Weir-type robust dual problem to the primal problem $(\hyperlink{CUP}{\mathrm{CUP}})$, and explore weak, strong, and converse duality properties. In addition, the obtained results are applied to an approximate uncertain multiobjective problem and a composite uncertain multiobjective problem with linear operators.
\end{abstract}
\textbf{Keywords}\hspace{3mm}Composite robust multiobjective optimization . Optimality conditions . Duality . Limiting subdifferential . Generalized convexity
\newline
\textbf{Mathematics Subject Classification (2020)}\hspace{3mm}49K99 . 65K10 . 90C29 . 90C46

\footnote{\textsf{$^*$Corresponding author}}

\afterpage{} \fancyhead{} \fancyfoot{} \fancyhead[LE, RO]{\bf\thepage} \fancyhead[LO]{\small Optimality and duality for composite uncertain multiobjective optimization in Asplund spaces} \fancyhead[RE]{\small Saadati and Oveisiha}
\section{Introduction}
\label{Sec1-Intro}
\textit{Composite optimization} problems have been widely recognized as a great important class of optimization problems. Their significance stems not only from their ability to cover a wide range of practical optimization problems, including minimax problems and penalty methods for constrained optimization, but also from their provision of a unified framework for studying the convergence behavior of various optimization algorithms and Lagrangian conditions (see, e.g., \citealt{Fletcher(1987)B, Burke(1985), Rockafellar(1988)}). \textit{Robust optimization} has become a powerful deterministic structure to study optimization problems under data uncertainty (\citealt{Beck(2009), Ben-Tal(2009)B, Ben-Tal(2002), Ben-Tal(2008), Jeyakumar(2015), Jeyakumar(2010), Jeyakumar(2012), Lee(2014)}). An \textit{uncertain optimization} problem usually associated with its robust counterpart which is known as the problem that the uncertain objective and constraint are satisfied for all possible scenarios within a prescribed uncertainty set. For classic contributions this field, we refer to \citet{Ben-Tal(2009)B}. Robust optimization approach considers the cases in which no probabilistic information about the uncertainties is given. In particular, most practical optimization problems often deal with uncertain data due to measurement errors, unforeseeable future developments, fluctuations, or disturbances, and depend on conflicting goals due to multiobjective decision makers which have different optimization criteria. So, the \textit{robust multiobjective optimization} is highly interesting in optimization theory and important in applications.

The first idea of robustness in multiobjective optimization was explored by \citet{Branke(1998)} and provided by \citet{Deb(2006)}. Here, the robustness concept is treated as a kind of sensitivity in the objective space against perturbations in the decision space. \citet{Kuroiwa(2012)} followed the robust approach (the worst-case approach) for multiobjective convex programming problems under uncertainty in both the objective functions and the constraints, and investigated necessary optimality conditions for weakly and properly robust efficient solutions. \citet{Ehrgott(2014)} extended the concept as presented by Kuroiwa and Lee, and interpreted a robust solution as a set of feasible solutions to the multiobjective problem of maximizing the objective function over the uncertainty set; see also the paper using the same approach (\citealt{Chuong(2020)}). After these works, \citet{Ide(2014)} derived various concepts of efficiency for uncertain multiobjective optimization problems, where only the objective functions were contaminated with different uncertain data, by replacing the set ordering with other set orderings, and presented numerical results on the occurrence of the various concepts.

Recently, \citet{Lee(2018)} dealt with robust multiobjective nonlinear semi-infinite programming with uncertain constraints, investigated necessary/sufficient conditions for weakly robust efficient solutions with the worst-case approach, and derived Wolfe-type dual problem and duality results. \citet{Chuong(2016)} considered uncertain multiobjective optimization problems involving nonsmooth/nonconvex functions, and introduced the concept of (strictly) generalized convexity  to establish optimality and duality theories with respect to limiting subdifferential for robust (weakly) Pareto solutions. \citet{Chen(2019)} studied necessary/sufficient conditions in terms of Clarke subdifferential for weakly and properly robust efficient solutions of nonsmooth multiobjective optimization problems with data uncertainty, formulated Mond-Weir-type dual problem and Wolfe-type dual problem, and explored duality results between the primal one and its dual problems under the generalized convexity assumptions. \citet{Fakhar(2018)} presented the nonsmooth sufficient optimality conditions for robust (weakly) efficient solutions and Mond-Weir-type duality results by applying the new concept of generalized convexity.

In addition to those stated above, the concept of approximate efficient solutions in multiobjective optimization problems, which can be viewed as feasible points whose objective values display a prescribed error $\varepsilon$ in the optimal values of the vector objective, has been studied widely. Optimality conditions and duality theories of $\varepsilon$-efficient solutions and $\varepsilon$-quasi-efficient solutions for various optimization problems under uncertainty have been presented by \citet{Chuong(2016)Positivity, Dutta(2005), Fakhar(2019), Govil(2004), Gupta(2008), Gutierrez(2006), Gutierrez(2010), Lee(2017), Lee(2012), Liu(1991), Saadati(2022)Optimization}. To the best of our knowledge, the most important results obtained in these directions did pay attention to finite-dimensional problems not dealing with composite functions. Hence, an infinite-dimensional setting is suitable to induce optimality and duality in composite optimization.

Suppose that $F : X \to W$ and $f : W \to Y$ be vector-valued functions between Asplund spaces, and that $K \subset Y$ be a pointed (i.e., $K \bigcap \, ( - K ) = \{ 0 \}$) closed convex cone with nonempty topological interior. Consider the following \textit{composite multiobjective optimization} problem of the form
\begin{equation*}
\hypertarget{CP}{}
  \begin{aligned}
    ( \mathrm{CP} ) \qquad \min\nolimits_{K} \,\,\, & ( f \circ F ) (x) \\
                           \mathrm{s.t.}     \,\,\, & ( g_{i} \circ G_{i} ) (x) \le 0, \quad i = 1, 2,\dots, n, \nonumber
  \end{aligned}
\end{equation*}
where $G = ( G_{1}, G_{2},\dots, G_{n} )$ and $g = ( g_{1}, g_{2},\dots, g_{n} )$ are vector-valued functions on Asplund spaces with components $G_{i} : X \to Z$ and $g_{i} : Z \to \mathbb{R}$, $i = 1, 2,\dots, n$. The problem $(\hyperlink{CP}{\mathrm{CP}})$ with data uncertainty in the constraints can be written by the \textit{composite uncertain multiobjective optimization} problem
\begin{equation*}
\hypertarget{CUP}{}
  \begin{aligned}
    ( \mathrm{CUP} ) \qquad \min\nolimits_{K} \,\,\, & ( f \circ F ) (x) \\
                            \mathrm{s.t.}     \,\,\, & ( g_{i} \circ G_{i} ) ( x, v_{i} ) \le 0, \quad i = 1, 2,\dots, n, \nonumber
  \end{aligned}
\end{equation*}
where $x \in X$ is a \textit{decision} variable, $v_{i}$ is an \textit{uncertain} parameter which belongs to the \textit{sequentially compact topological} space $\mathcal{V}_{i}$, and $G_{i} : X \times \mathcal{V}_{i} \to Z \times \mathcal{U}_{i}$ and $g_{i} : Z \times \mathcal{U}_{i} \to \mathbb{R}$ are given functions for topological space $\mathcal{U}_{i}$, $i = 1, 2,\dots, n$.

For investigating the problem $(\hyperlink{CUP}{\mathrm{CUP}})$, we associate with it the so-called \textit{robust} counterpart
\begin{equation*}
\hypertarget{CRP}{}
  \begin{aligned}
    ( \mathrm{CRP} ) \qquad \min\nolimits_{K} \,\,\, & ( f \circ F ) (x) \\
                            \mathrm{s.t.}     \,\,\, & ( g_{i} \circ G_{i} ) ( x, v_{i} ) \le 0, \quad \forall v_{i} \in \mathcal{V}_{i}, \,\, i = 1, 2,\dots, n.
  \end{aligned}
\end{equation*}
Let
\begin{equation*}
   C:= \medBl x \in X \,\mid\, ( g_{i} \circ G_{i} ) ( x, v_{i} ) \le 0, \,\, \forall v_{i} \in \mathcal{V}_{i}, \,\, i = 1, 2,\dots, n \medBr
\end{equation*}
be the \textit{feasible} set of the problem $(\hyperlink{CRP}{\mathrm{CRP}})$.
\begin{rem}
\label{Rem1-1}
The problem $(\hyperlink{CUP}{\mathrm{CUP}})$ provides a quite general framework for various uncertain multiobjective optimization problems as follows
\begin{itemize}
  \item[(i)] if $X = W = Z$, $\mathcal{V}_{i} = \mathcal{U}_{i}$, $i = 1, 2,\dots, n$, and $F$ and $G$ are identical maps, the problem $(\hyperlink{CUP}{\mathrm{CUP}})$ collapses to the following \textit{uncertain multiobjective optimization} problem stated by \citet{Saadati(2022)}
      \begin{equation*}
      \hypertarget{UP}{}
        ( \mathrm{UP} ) \qquad \min\nolimits_{K} \,\,\, \medBl f (x) \,\mid\, x \in X, \,\, g_{i} ( x, v_{i} ) \le 0, \,\, i = 1, 2,\dots, n \medBr.
      \end{equation*}
  \item[(ii)] if $X = W = Z$, $\mathcal{V}_{i} = \mathcal{U}_{i}$, $i = 1, 2,\dots, n$, $\mathcal{V}_{n+j} = \mathcal{U}_{n+j}$, $j = 1, 2,\dots, m$, $F$ and $G$ are identical maps, and $g$ is given by
      \begin{align}
        g ( x, v ) := ( g_{1} ( x, v_{1} ), g_{2} ( x, v_{2} ),\dots, g_{n} ( x, v_{n} )&, h_{1} ( x, v_{n+1} ), h_{2} ( x, v_{n+2} ),\dots, h_{m} ( x, v_{n+m} ) ), \nonumber \\
        &x \in X, \nonumber \\
        &v_{i} \in \mathcal{V}_{i}, \,\, i = 1, 2,\dots, n, \nonumber \\
        \label{1-1}
        &v_{n+j} \in \mathcal{V}_{n+j}, \,\, j = 1, 2,\dots, m,
      \end{align}
      then the problem $(\hyperlink{CUP}{\mathrm{CUP}})$ reduces to a (\textit{standard}) \textit{uncertain multiobjective optimization} problem of the form
      \begin{equation*}
      \hypertarget{SUP}{}
        \begin{aligned}
          ( \mathrm{SUP} ) \qquad \min\nolimits_{K} \,\,\, \medBl f (x) \,\mid\, x \in X, \,\, &g_{i} ( x, v_{i} ) \le 0, \,\, i = 1, 2,\dots, n, \\
          &h_{j} ( x, v_{n+j} ) = 0, \,\, j = 1, 2,\dots, m \medBr.
        \end{aligned}
      \end{equation*}
  \item[(iii)] if $X = W = Z$, $Y := \mathbb{R}^{p}$, $K := \mathbb{R}^{p}_{+}$, $\mathcal{V}_{i} = \mathcal{U}_{i}$, $ i = 1, 2,\dots, m$, and $F$ and $G$ are identical maps, the problem $(\hyperlink{CUP}{\mathrm{CUP}})$ collapses to an uncertain multiobjective optimization problem defined by \citet{Fakhar(2018)}.
  \item[(iv)] if $X = W = Z := \mathbb{R}^{n}$, $Y := \mathbb{R}^{m}$, $K := \mathbb{R}^{m}_{+}$, $\mathcal{V}_{i} = \mathcal{U}_{i}$, $i = 1, 2,\dots, l$, and $F$ and $G$ are identical maps, then the problem $(\hyperlink{CUP}{\mathrm{CUP}})$ reduces to an uncertain multiobjective optimization problem presented by \citet{Chuong(2016)}.
\end{itemize}
\end{rem}
\begin{defn}
\hypertarget{Def1-1}{}
\begin{itemize}
  \item[(i)] We say that a vector $\bar{x} \in X$ is a \textit{robust efficient solution} of the problem $(\hyperlink{CUP}{\mathrm{CUP}})$, and write $\bar{x} \in \mathcal{S}(\hyperlink{CRP}{\mathrm{CRP}})$, if $\bar{x}$ is an \textit{efficient solution} of the problem (\hyperlink{CRP}{CRP}), i.e., $\bar{x} \in C$ and
      \begin{equation*}
        ( f \circ F ) (x) - ( f \circ F ) ( \bar{x} ) \notin - K \setminus \{ 0 \}, \quad \forall x \in C.
      \end{equation*}
  \item[(ii)] A vector $\bar{x} \in X$ is called a \textit{weakly robust efficient solution} of the problem $(\hyperlink{CUP}{\mathrm{CUP}})$, and write $\bar{x} \in \mathcal{S}^{w}(\hyperlink{CRP}{\mathrm{CRP}})$, if $\bar{x}$ is a \textit{weakly efficient solution} of the problem (\hyperlink{CRP}{CRP}), i.e., $\bar{x} \in C$ and
      \begin{equation*}
        ( f \circ F ) (x) - ( f \circ F ) ( \bar{x} ) \notin - \mathrm{int}\hspace{.4mm}K, \quad \forall x \in C.
      \end{equation*}
\end{itemize}
\end{defn}

The rest of this paper is organized as follows. Section \ref{Sec2-Preli} contains some preliminary definitions and several auxiliary results from nonsmooth variational analysis. In Section \ref{Sec3-NecSuf}, we establish necessary/sufficient optimality conditions for the existence of weakly robust efficient solutions and also sufficient conditions for robust efficient solutions of the problem $(\hyperlink{CUP}{\mathrm{CUP}})$ in terms of the limiting subdifferential. Section \ref{Sec4-Dual} is concerned with the duality relations for (weakly) robust efficient solutions between the corresponding problems. The concluding Section \ref{Sec5-Appl} provides applications of special composite forms to the robust multiobjective optimization.
\section{Preliminaries}
\label{Sec2-Preli}
Throughout this paper, we use the standard notation of variational analysis; see, for example, \citealt{Mordukhovich(2006)B}. Unless otherwise stated, all the spaces under consideration are \textit{Asplund} with the norm $\| \cdot \|$ and the canonical pairing $\langle \cdot \,, \cdot \rangle$ between the space $X$ in question and its \textit{dual} $X^{*}$ equipped with the \textit{weak}$^{*}$ \textit{topology} $w^{*}$. By $B_{X} ( x, r )$, we denote the \textit{closed ball} centered at $x \in X$ with radius $r > 0$, while $B_{X}$ and $B_{ X^{*} }$ stand for the \textit{closed unit ball} in $X$ and $X^{*}$, respectively. For a given nonempty set $\Omega \subset X$, the symbols $\mathrm{co}\hspace{.4mm}\Omega$, $\mathrm{cl}\hspace{.4mm}\Omega$, and $\mathrm{int}\hspace{.4mm}\Omega$ indicate the \textit{convex hull}, \textit{topological closure}, and \textit{topological interior} of $\Omega$, respectively, while $\mathrm{cl}^{*}\Omega$ stands for the \textit{weak}$^{*}$ \textit{topological closure} of $\Omega \subset X^{*}$. The \textit{dual cone} of $\Omega$ is the set
\begin{equation*}
  \Omega^{+} := \medBl x^{*} \in X^{*} \,\mid\, \langle x^{*}, x \rangle \ge 0, \,\,\, \forall x \in \Omega \medBr.
\end{equation*}
For $n \in \mathbb{N} := \{ 1, 2,\dots \}$, $\mathbb{R}^{n}_{+}$ denotes the nonnegative orthant of $\mathbb{R}^{n}$. Besides, the symbol $T^{\top}$ signifies the adjoint operator or conjugate transpose of the linear operator $T$.

A given set-valued mapping $H : \Omega \subset X \overrightarrow{\to} X^{*}$ is called \textit{weak}$^{*}$ \textit{closed} at $\bar{x} \in \Omega$ if for any sequence $\{ x_{k} \} \subset \Omega$, $x_{k} \to \bar{x}$, and any sequence $\{ x^{*}_{k} \} \subset X^{*}$, $x^{*}_{k} \in H ( x_{k} )$, $x^{*}_{k} \overset{w^{*}} \to x^{*}$, one has $x^{*} \in H ( \bar{x} )$.

For a set-valued mapping $H : X \overrightarrow{\to} X^{*}$, the \textit{sequential Painlev\'{e}-Kuratowski upper/outer limit} of $H$ as $x \to \bar{x}$ is defined by
\begin{align*}
  \underset{x \to \bar{x}}{\mathrm{Lim}\sup} \, H (x) := \medBl x^{*} \in X^{*} \,\mid\,\,\, &\exists \text{ sequences } x_{k} \to \bar{x} \text{ and } x^{*}_{k} \overset{\scriptscriptstyle w^{*}} \to x^{*} \\
  &\text{with } x^{*}_{k} \in H ( x_{k} ) \text{ for all } k \in \mathbb{N} \medBr.
\end{align*}

Let $\Omega \subset X$ be \textit{locally closed} around $\bar{x} \in \Omega$, i.e., there is a neighborhood $U$ of $\bar{x}$ for which $\Omega \bigcap \mathrm{cl}\hspace{.4mm}U$ is closed. The \textit{Fr\'{e}chet normal cone} $\widehat{N} ( \bar{x}; \Omega )$ and the \textit{Mordukhovich normal cone} $N ( \bar{x}; \Omega )$ to $\Omega$ at $\bar{x} \in \Omega$ are defined, respectively, by
\begin{align}
\label{2-1}
  \widehat{N} ( \bar{x}; \Omega ) &:= \medBl x^{*} \in X^{*} \,\mid\, \limsup\limits_{ x \overset{\hspace{-1mm}\sOmega} \to \bar{x} } \dfrac{ \langle x^{*}, x - \bar{x} \rangle }{ \| x - \bar{x} \| } \le 0 \medBr \\
\intertext{and}
\label{2-2}
  N ( \bar{x}; \Omega ) &:= \underset{x \overset{\hspace{-1mm}\sOmega} \to \bar{x}}{\mathrm{Lim}\sup} \, \widehat{N} ( x; \Omega ),
\end{align}
where $x \overset{\hspace{-1mm}\Omega} \to \bar{x}$ stands for $x \to \bar{x}$ with $x \in \Omega$. If $\bar{x} \notin \Omega$, we put $\widehat{N} ( \bar{x}; \Omega ) = N ( \bar{x}; \Omega ) := \emptyset$.

For an extended real-valued function $\phi : X \to \overline{\mathbb{R}}$, the \textit{limiting/Mordukhovich subdifferential} and the \textit{regular/Fr\'{e}chet subdifferential} of $\phi$ at $\bar{x} \in \mathrm{dom}\,\phi$ are given, respectively, by
\begin{align*}
  \partial \phi ( \bar{x} ) &:= \medBl x^{*} \in X^{*} \,\mid\, ( x^{*}, - 1 ) \in  N ( ( \bar{x}, \phi (x) ); \mathrm{epi}\,\phi ) \medBr \\
\intertext{and}
  \widehat{\partial} \phi ( \bar{x} ) &:= \medBl x^{*} \in X^{*} \,\mid\, ( x^{*}, - 1 ) \in  \widehat{N} ( ( \bar{x}, \phi (x) ); \mathrm{epi}\,\phi ) \medBr.
\end{align*}
If $| \phi ( \bar{x} ) | = \infty$, then one puts $\partial \phi ( \bar{x} ) = \widehat{\partial} \phi ( \bar{x} ) := \emptyset$.

For a vector-valued function $f : X \to Y$, we apply a scalarization formula with respect to some $y^{*} \in Y^{*}$ defined by
\begin{equation*}
  \langle y^{*}, f \rangle (x) := \langle y^{*}, f (x) \rangle, \quad x \in X.
\end{equation*}
\noindent
We recall another expression of the scalarization scheme in the next lemma.
\begin{lem}
\label{Lem2-1}
Let $y^{*} \in Y^{*}$, and let $f : X \to Y$ be Lipschitz continuous around $\bar{x} \in X$. We have
\begin{itemize}
  \item[{\normalfont (i)}] {\normalfont (See (\citealt[Proposition~3.5]{Mordukhovich(2006)}))} $x^{*} \in \widehat{\partial} \langle y^{*}, f \rangle ( \bar{x} ) \,\, \Leftrightarrow \,\, ( x^{*}, - y^{*} ) \in \widehat{N} ( ( \bar{x}, f ( \bar{x} ) ); \mathrm{gph}\,f )$.
  \item[{\normalfont (ii)}] {\normalfont (See (\citealt[Theorem~1.90]{Mordukhovich(2006)B}))} $x^{*} \in \partial \langle y^{*}, f \rangle ( \bar{x} ) \,\, \Leftrightarrow \,\, ( x^{*}, - y^{*} ) \in N ( ( \bar{x}, f ( \bar{x} ) ); \mathrm{gph}\,f )$.
\end{itemize}
\end{lem}
\noindent
The following lemma gives a \textit{chain rule} for the limiting subdifferential.
\begin{lem}
\label{Lem2-2}
{\normalfont (See (\citealt[Corollary~3.43]{Mordukhovich(2006)B}))}
Let $f : X \to Y$ be locally Lipschitz at $\bar{x} \in X$, and let $\phi : Y \to \mathbb{R}$ be locally Lipschiz at $f ( \bar{x} )$. Then one has
\begin{equation*}
  \partial ( \phi \circ f ) ( \bar{x} ) \subset \medcup\limits_{ y^{*} \in \partial \phi ( f ( \bar{x} ) ) } \partial \langle y^{*}, f \rangle ( \bar{x} ).
\end{equation*}
\end{lem}
\noindent
The \textit{sum rule} for the limiting subdifferential will be useful in our analysis.
\begin{lem}
\label{Lem2-3}
{\normalfont (See (\citealt[Theorem~3.36]{Mordukhovich(2006)B}))}
Let $\phi_{i} : X \to \overline{\mathbb{R}}$, $( i \in \{ 1, 2,\dots, n \}, n \ge 2 )$, be lower semicontinuous around $\bar{x}$, and let all but one of these functions be Lipschitz continuous around $\bar{x} \in X$. Then one has
\begin{equation*}
  \partial ( \phi_{1} + \phi_{2} +\dots+ \phi_{n} ) ( \bar{x} ) \subset \partial \phi_{1} ( \bar{x} ) + \partial \phi_{2} ( \bar{x} ) +\dots+ \partial \phi_{n} ( \bar{x} ).
\end{equation*}
\end{lem}

It is worth to mention that inspecting the proof of (\citealt[Theorem~3.3]{Chuong(2016)}) (see also \citealt{Mordukhovich(2013)MathOperRes, Mordukhovich(2013)SIAMJOptim}) reveals that this proof contains a formula for the limiting subdifferential of \textit{maximum} functions in finite-dimensional spaces. The following lemma generalizes the corresponding result in arbitrary Asplund spaces. Its proof is on the straightforward side and similar in some places given in \citealt{Chuong(2016)}, and so we omit the details. The notation $\partial_{x}$ signifies the limiting subdifferential operation with respect to $x$.
\begin{lem}
\label{Lem2-4}
Let $\mathcal{V}$ be a sequentially compact topological space, and let $g : X \times \mathcal{V} \to \mathbb{R}$ be a function such that for each fixed $v \in \mathcal{V}$, $g ( \cdot, v )$ is Lipschitz continuous around $\bar{x} \in X$ and $g ( \bar{x}, \cdot )$ is upper semicontinuous on $\mathcal{V}$. Let $\phi (x) := \max\limits_{ v \in \mathcal{V} } g ( x, v )$. If the multifunction $( x, v ) \in X \times \mathcal{V} \,\, \overrightarrow{\to} \,\, \partial_{x} g ( x, v ) \subset X^{*}$ is weak$^{*}$ closed at $( \bar{x}, \bar{v} )$ for each $\bar{v} \in \mathcal{V} ( \bar{x} )$, then the set $\mathrm{cl}^{*}\mathrm{co} \medPl \medcup \medBl \partial_{x} g ( \bar{x}, v ) \,\mid\, v \in \mathcal{V} ( \bar{x} ) \medBr \medPr$ is nonempty and
\begin{equation*}
  \partial \phi ( \bar{x} ) \subset \mathrm{cl}^{*}\mathrm{co} \medPl \medcup \medBl \partial_{x} g ( \bar{x}, v ) \,\mid\, v \in \mathcal{V} ( \bar{x} ) \medBr \medPr,
\end{equation*}
where $\mathcal{V} ( \bar{x} ) := \medBl v \in \mathcal{V} \,\mid\, g ( \bar{x}, v ) = \phi ( \bar{x} ) \medBr$.
\end{lem}
\noindent
The next lemma is concerning with the limiting subdifferential for the maximum of finitely many functions in Asplund spaces.
\begin{lem}
\label{Lem2-5}
{\normalfont (See (\citealt[Theorem~3.46]{Mordukhovich(2006)B}))}
Let $\phi_{i} : X \to \overline{\mathbb{R}}$, $( i \in \{ 1, 2,\dots, n \}, n \ge 2 )$, be Lipschitz continuous around $\bar{x}$. Put $\phi (x) := \max\limits_{ i \in \{ 1, 2,\dots, n \} } \phi_{i} (x)$. Then
\begin{equation*}
  \partial \phi ( \bar{x} ) \subset \medcup \medBl \partial \medPl \sum_{ i \in I ( \bar{x} ) } \mu_{i} \, \phi_{i} \medPr ( \bar{x} ) \,\mid\, ( \mu_{1}, \mu_{2},\dots, \mu_{n} ) \in \Lambda ( \bar{x} ) \medBr,
\end{equation*}
where
\begin{equation*}
  I ( \bar{x} ) := \medBl i \in \{ 1, 2,\dots, n \} \,\mid\, \phi_{i} ( \bar{x} ) = \phi ( \bar{x} ) \medBr
\end{equation*}
and
\begin{equation*}
  \Lambda ( \bar{x} ) := \medBl ( \mu_{1}, \mu_{2},\dots, \mu_{n} ) \,\mid\, \mu_{i} \ge 0, \,\, \sum_{i=1}^{n} \mu_{i} = 1, \,\, \mu_{i} \, ( \phi_{i} ( \bar{x} ) - \phi ( \bar{x} ) ) = 0 \medBr.
\end{equation*}
\end{lem}
\noindent
The following lemma computes the limiting subdifferential of a norm.
\begin{lem}
\label{Lem2-6}
{\normalfont (See (\citealt[Lemma~4.1.11]{Gopfert(2003)B1}))}
Let $x \in X$ and $\beta > 1$. Then we have
\begin{align*}
  &\partial \| x \| =
  \begin{cases}
    \medBl x^{*} \in X^{*} \,\mid\, \langle x^{*}, x \rangle = \| x \|, \,\, \| x^{*} \| = 1 \medBr & \text{if } x \ne 0, \\
    \medBl x^{*} \in X^{*} \,\mid\, \| x^{*} \| \le 1 \medBr & \text{if } x = 0,
  \end{cases}\\
  \intertext{and}
  &\partial \bigg( \dfrac{1}{\beta} \| x \|^{\beta} \bigg) = \medBl x^{*} \in X^{*} \,\mid\, \langle x^{*}, x \rangle = \| x \|^{\beta}, \,\, \| x^{*} \| = \| x \|^{\beta - 1} \medBr.
\end{align*}
\end{lem}
\begin{assum}
\hypertarget{Assu2-1}{}
{\normalfont (See (\citealt[p.131]{Chuong(2016)}))}
\begin{itemize}
  \item[(A1)] For a fixed $\bar{x} \in X$, $F$ is locally Lipschitz at $\bar{x}$ and $f$ is locally Lipschitz at $F ( \bar{x} )$.
  \item[(A2)] For each $i = 1, 2, \dots, n$, $G_{i}$ is locally Lipschitz at $\bar{x}$ and uniformly on $\mathcal{V}_{i}$, and $g_{i}$ is Lipschitz continuous on $G_{i} ( \bar{x}, \mathcal{V}_{i} )$.
  \item[(A3)] For each $i = 1, 2, \dots, n$, the functions $v_{i} \in \mathcal{V}_{i} \mapsto G_{i} ( \bar{x}, v_{i} ) \in Z \times \mathcal{U}_{i}$ and $G_{i} ( \bar{x}, v_{i} ) \mapsto g_{i} ( G_{i} ( \bar{x}, v_{i} ) ) \in \mathbb{R}$ are locally Lipschitzian.
  \item[(A4)] For each $i = 1, 2, \dots, n$, we define the real-valued functions $\phi_{i}$ and $\phi$ on $X$ via
      \begin{equation*}
        \phi_{i} (x) := \max_{ v_{i} \in \mathcal{V}_{i} } ( g_{i} \circ G_{i} ) ( x, v_{i} ) \qquad \text{and} \qquad \phi (x) := \max_{ i \in \{ 1, 2,\dots, n \} } \phi_{i} (x),
      \end{equation*}
      and we notice that the assumption (A3) implies that $\phi_{i}$ is well defined on $\mathcal{V}_{i}$. In addition, $\phi_{i}$ and $\phi$ follow readily that are locally Lipschitz at $\bar{x}$, since each $( g_{i} \circ G_{i} ) ( \bar{x}, v_{i} )$ is (see, e.g., (\citealt[(H1), p.131]{Chuong(2016)}) or (\citealt[p.290]{Lee(2014)})). Note that the feasible set $C$ can be equivalently characterized by
      \begin{equation*}
        C = \medBl x \in X \,\mid\, \phi_{i} (x) \le 0, \,\, i = 1, 2, \dots, n \medBr = \medBl x \in X \,\mid\, \phi (x) \le 0 \medBr.
      \end{equation*}
  \item[(A5)] For each $i = 1, 2, \dots, n$, the multifunction $( x, v_{i} ) \in X \times \mathcal{V}_{i} \,\, \overrightarrow{\to} \,\, \partial_{x} ( g_{i} \circ G_{i} ) ( x, v_{i} ) \subset X^{*}$ is weak$^{*}$ closed at $( \bar{x}, \bar{v}_{i} )$ for each $\bar{v}_{i} \in \mathcal{V}_{i} ( \bar{x} )$, where $\mathcal{V}_{i} ( \bar{x} ) = \medBl v_{i} \in \mathcal{V}_{i} \,\mid\, ( g_{i} \circ G_{i} ) ( \bar{x}, v_{i} ) = \phi_{i} ( \bar{x} ) \medBr$.
\end{itemize}
\end{assum}
\section{Robust necessary and sufficient optimality}
\label{Sec3-NecSuf}
In this section, we study optimality conditions in composite robust multiobjective optimization problems. More precisely, first by exploiting the nonsmooth version of Fermat's rule, sum rule, and chain rule for the limiting subdifferential, necessary conditions for weakly robust efficient solutions of the problem $(\hyperlink{CUP}{\mathrm{CUP}})$ will be established. We then derive sufficient conditions for the existence of such solutions as well as robust efficient solutions under assumptions of pseudo-quasi convexity for composite vector-valued functions.

The first theorem in this section provides a necessary optimality condition in the sense of the limiting subdifferential for weakly robust efficient solutions of the problem $(\hyperlink{CUP}{\mathrm{CUP}})$. To prove, we need to state a \textit{fuzzy} necessary condition expressed in terms of the Fr\'{e}chet subdifferential for weakly robust efficient solutions of the problem $(\hyperlink{UP}{\mathrm{UP}})$ as follows.
\begin{thm}
\label{Thm3-1}
{\normalfont (See (\citealt[Theorem~3.1]{Chuong(2019)}))}
Let $\bar{x}$ be a weakly robust efficient solution of the problem $(\hyperlink{UP}{\mathrm{UP}})$. Then for each $k \in \mathbb{N}$ there exist $x^{1k} \in B_{X} ( \bar{x}, \frac{1}{k} )$, $x^{2k} \in B_{X} ( \bar{x}, \frac{1}{k} )$, $y^{*}_{k} \in K^{+}$ with $\| y^{*}_{k} \| = 1$, and $\alpha_{k} \in \mathbb{R}_{+}$ such that
\begin{align*}
  & 0 \in \widehat{\partial} \langle y_{k}^{*}, f \rangle ( x^{1k} ) + \alpha_{k} \, \widehat{\partial} \phi ( x^{2k} ) + \dfrac{1}{k} B_{X^{*}}, \nonumber \\
  & | \alpha_{k} \, \phi ( x^{2k} ) | \le \dfrac{1}{k}.
\end{align*}
\end{thm}
\begin{thm}
\label{Thm3-2}
Suppose that $\bar{x} \in \mathcal{S}^{w}(\hyperlink{CRP}{\mathrm{CRP}})$. Then there exist $y^{*} \in K^{+}$, $\mu := ( \mu_{1}, \mu_{2},\dots, \mu_{n} ) \in \mathbb{R}^{n}_{+}$, with $ \| y^{*} \| + \| \mu \| = 1$, and $\bar{v}_{i} \in \mathcal{V}_{i}$, $i = 1, 2,\dots, n$, such that
\begin{equation}
\label{3-1}
  \left\{
    \begin{aligned}
      & 0 \in \medcup\limits_{ w^{*} \in \partial \langle y^{*}, f \rangle ( F ( \bar{x} ) ) } \partial \langle w^{*}, F \rangle ( \bar{x} ) + \sum_{i=1}^{n} \mu_{i} \, \mathrm{cl}^{*}\mathrm{co} \medPl \medcup \medBl \medcup\limits_{ v_{i}^{*} \in \partial_{x} g_{i} ( G_{i} ( \bar{x}, v_{i} ) ) } \partial_{x} \langle v_{i}^{*}, G_{i} \rangle ( \bar{x}, v_{i} ) \,\mid\, v_{i} \in \mathcal{V}_{i} ( \bar{x} ) \medBr \medPr, \\
      & \mu_{i} \, \max\limits_{ v_{i} \in \mathcal{V}_{i} } g_{i} ( G_{i} ( \bar{x}, v_{i} ) ) = \mu_{i} \, g_{i} ( G_{i} ( \bar{x}, \bar{v}_{i} ) ) = 0, \quad i = 1, 2,\dots, n.
    \end{aligned}
  \right.
\end{equation}
\end{thm}
\begin{proof}
Let us put $\tilde{f} := f \circ F$. In this case, the problem $(\hyperlink{CRP}{\mathrm{CRP}})$ changes to
\begin{equation*}
\hypertarget{CRPt}{}
  \begin{aligned}
    ( \overset{\sim}{\mathrm{CRP}} ) \qquad \min\nolimits_{K} \,\,\, \medBl \tilde{f} (x) \,\mid\, x \in X, \,\, \phi (x) \le 0 \medBr.
  \end{aligned}
\end{equation*}
Using Theorem \ref{Thm3-1} to the problem $(\hyperlink{CRPt}{\overset{\sim}{\mathrm{CRP}}})$, there exist sequences $x^{1k} \to \bar{x}$, $x^{2k} \to \bar{x}$, $y^{*}_{k} \in K^{+}$ with $\| y^{*}_{k} \| = 1$, $\alpha_{k} \in \mathbb{R}_{+}$, $x^{*}_{1k} \in \widehat{\partial} \langle y^{*}_{k}, f \circ F \rangle ( x^{1k} )$, and $x^{*}_{2k} \in \alpha_{k} \, \widehat{\partial} \phi ( x^{2k} )$ satisfying
\begin{align}
\label{3-2}
  & 0 \in x^{*}_{1k} + x^{*}_{2k} + \dfrac{1}{k} B_{ X^{*} }, \\
  & \alpha_{k} \, \phi ( x^{2k} ) \to 0 \text{ as } k \to \infty. \nonumber
\end{align}
We now consider the following two possibilities

\textbf{Case 1:} If $\{ \alpha_{k} \}$ is bounded, then without loss of generality we may suppose that $\alpha_{k} \to \alpha \in \mathbb{R}_{+}$ as $k \to \infty$. Besides, as the sequence $\{ y^{*}_{k} \} \subset K^{+}$ is bounded, by invoking the weak$^{*}$ sequential compactness of bounded sets in duals to Asplund spaces, we don't restrict the generality by assuming that $y^{*}_{k} \overset{w^{*}} \rightarrow \bar{y}^{*} \in K^{+}$ with $\| \bar{y}^{*} \| = 1$ as $k \to \infty$. Let $\ell > 0$ be a constant modulus for the locally Lipschitz function $f \circ F$ at $\bar{x}$. It is clear that $\| x_{1k}^{*} \| \le \ell \, \| y^{*}_{k} \| \le \ell$ for all $k \in \mathbb{N}$ (see (\citealt[Proposition~1.85]{Mordukhovich(2006)B})). As above, by passing to a subsequence if necessary, that $x^{*}_{1k} \overset{w^{*}} \rightarrow x^{*}_{1} \in X^{*}$ as $k \to \infty$, and so it follows from (\ref{3-2}) that $x^{*}_{2k} \overset{w^{*}} \rightarrow x^{*}_{2} := - x^{*}_{1}$ as $k \to \infty$. According to Lemma \ref{Lem2-1}, we deduce from the inclusion $x^{*}_{1k} \in \widehat{\partial} \langle y^{*}_{k}, f \circ F \rangle ( x^{1k} )$ that
\begin{equation*}
  ( x_{1k}^{*}, - y_{k}^{*} ) \in \widehat{N} ( ( x^{1k}, ( f \circ F ) ( x^{1k} ) ); \mathrm{gph}\,( f \circ F ) ), \quad k \in \mathbb{N}.
\end{equation*}
Letting $k \to \infty$ and noticing the definitions (\ref{2-1}) and (\ref{2-2}) of the normal cones, we obtain the relation $( x^{*}_{1}, -\bar{y}^{*} ) \in N ( ( \bar{x}, ( f \circ F ) ( \bar{x} ) ); \mathrm{gph}\,( f \circ F ) )$ which is equivalent to
\begin{equation}
\label{3-3}
  x^{*}_{1} \in \partial \langle \bar{y}^{*}, f \circ F \rangle ( \bar{x} ),
\end{equation}
because of Lemma \ref{Lem2-1}. Similarly, we get $x^{*}_{2} \in \alpha \, \partial \phi ( \bar{x} )$. This combined with (\ref{3-3}) entails that
\begin{equation}
\label{3-4}
  0 \in \partial \langle \bar{y}^{*}, f \circ F \rangle ( \bar{x} ) + \alpha \, \partial \phi ( \bar{x} ),
\end{equation}
by taking $x^{*}_{2} = - x^{*}_{1}$. Applying Lemma \ref{Lem2-5} to the limiting subdifferential of the maximum function $\phi$, we arrive at
\begin{equation}
\label{3-5}
  \partial \phi ( \bar{x} ) \subset \medcup \medBl \partial \medPl \sum_{ i \in I ( \bar{x} ) } \mu_{i} \, \phi_{i} \medPr ( \bar{x} ) \,\mid\, ( \mu_{1}, \mu_{2},\dots, \mu_{n} ) \in \Lambda ( \bar{x} ) \medBr,
\end{equation}
where $I ( \bar{x} ) = \medBl i \in \{ 1, 2,\dots, n \} \,\mid\, \phi_{i} ( \bar{x} ) = \phi ( \bar{x} ) \medBr$ and
\begin{equation*}
  \Lambda ( \bar{x} ) = \medBl ( \mu_{1}, \mu_{2},\dots, \mu_{n} ) \,\mid\, \mu_{i} \ge 0, \,\, \sum_{i=1}^{n} \mu_{i} = 1, \,\, \mu_{i} \, ( \phi_{i} ( \bar{x} ) - \phi ( \bar{x} ) ) = 0 \medBr.
\end{equation*}
Invoking further Lemma \ref{Lem2-4} allows us
\begin{equation}
\label{3-6}
  \partial \phi_{i} ( \bar{x} ) \subset \mathrm{cl}^{*}\mathrm{co} \medPl \medcup \medBl \partial_{x} ( g_{i} \circ G_{i} ) ( \bar{x}, v_{i} ) \,\mid\, v_{i} \in \mathcal{V}_{i} ( \bar{x} ) \medBr \medPr, \quad i = 1, 2,\dots, n,
\end{equation}
where $\mathcal{V}_{i} ( \bar{x} ) = \medBl v_{i} \in \mathcal{V}_{i} \,\mid\, ( g_{i} \circ G_{i} ) ( \bar{x}, v_{i} ) = \phi_{i} ( \bar{x} ) \medBr$ and the set $\mathrm{cl}^{*}\mathrm{co} \medPl \medcup \medBl \partial_{x} ( g_{i} \circ G_{i} ) ( \bar{x}, v_{i} ) \,\mid\, v_{i} \in \mathcal{V}_{i} ( \bar{x} ) \medBr \medPr$ is nonempty. By the sum rule of Lemma \ref{Lem2-3}, it follows from the relations (\ref{3-4})-(\ref{3-6}) that
\begin{align*}
  0 &\in \partial \langle \bar{y}^{*}, f \circ F \rangle ( \bar{x} ) + \alpha \, \medcup \medBl \sum_{ i \in I ( \bar{x} ) } \mu_{i} \, \mathrm{cl}^{*}\mathrm{co} \medPl \medcup \medBl \partial_{x} ( g_{i} \circ G_{i} ) ( \bar{x}, v_{i} ) \,\mid\, v_{i} \in \mathcal{V}_{i} ( \bar{x} ) \medBr \medPr \,\mid\, \\
  &( \mu_{1}, \mu_{2},\dots, \mu_{n} ) \in \Lambda ( \bar{x} ) \medBr.
\end{align*}
Thus, there exist $\bar{ \mu } := ( \bar{\mu}_{1}, \bar{\mu}_{2},\dots, \bar{\mu}_{n} ) \in \Lambda ( \bar{x} )$, with $\mathlarger{\sum\limits }_{i=1}^{n} \, \bar{\mu}_{i} = 1$ and $\bar{\mu}_{i} = 0$ for all $i \in \{ 1, 2,\dots, n \} \setminus I ( \bar{x} )$, satisfying
\begin{equation*}
  0 \in \partial \langle \bar{y}^{*}, f \circ F \rangle ( \bar{x} ) + \alpha \, \sum_{i=1}^{n} \bar{\mu}_{i} \, \mathrm{cl}^{*}\mathrm{co} \medPl \medcup \medBl \partial_{x} ( g_{i} \circ G_{i} ) ( \bar{x}, v_{i} ) \,\mid\, v_{i} \in \mathcal{V}_{i} ( \bar{x} ) \medBr \medPr.
\end{equation*}
Dividing the above inclusion by $c := \| \bar{y}^{*} \| + \alpha \, \| \bar{\mu} \|$ and then letting $y^{*} := \dfrac{ \bar{y}^{*} }{c}$ and $\mu := \dfrac{\alpha}{c} \, \bar{\mu}$, we have $y^{*} \in K^{+}$ and $\mu := ( \mu_{1}, \mu_{2},\dots, \mu_{n} ) \in \mathbb{R}^{n}_{+}$, with $\| y^{*} \| + \| \mu \| = 1 $, such that
\begin{equation*}
  0 \in \partial \langle y^{*}, f \circ F \rangle ( \bar{x} ) + \sum_{i=1}^{n} \mu_{i} \, \mathrm{cl}^{*}\mathrm{co} \medPl \medcup \medBl \partial_{x} ( g_{i} \circ G_{i} ) ( \bar{x}, v_{i} ) \,\mid\, v_{i} \in \mathcal{V}_{i} ( \bar{x} ) \medBr \medPr.
\end{equation*}
Putting now $\psi := \langle y^{*}, f \rangle$, we rewrite the above inclusion as
\begin{equation}
\label{3-7}
  0 \in \partial ( \psi \circ F ) ( \bar{x} ) + \sum_{i=1}^{n} \mu_{i} \, \mathrm{cl}^{*}\mathrm{co} \medPl \medcup \medBl \partial_{x} ( g_{i} \circ G_{i} ) ( \bar{x}, v_{i} ) \,\mid\, v_{i} \in \mathcal{V}_{i} ( \bar{x} ) \medBr \medPr.
\end{equation}
By the assumptions (\hyperlink{Assu2-1}{A1}) and (\hyperlink{Assu2-1}{A2}), $F$ and $\psi$ are locally Lipschitzian at $\bar{x}$ and $F ( \bar{x} )$, respectively. And $G_{i}$ is locally Lipschitzian at $\bar{x}$ and uniformly on $\mathcal{V}_{i}$, and $g_{i}$ is Lipschitz continuous on $G_{i} ( \bar{x}, \mathcal{V}_{i} )$ for each $i = 1, 2,\dots, n$. So using the limiting subdifferential chain rule of Lemma \ref{Lem2-2} for (\ref{3-7}), we get the first relation in this theorem.

On the other side, due to the locally Lipschitz continuity of the function $v_{i} \in \mathcal{V}_{i} \longmapsto ( g_{i} \circ G_{i} ) ( \bar{x}, v_{i} )$ and the sequentially compactness of $\mathcal{V}_{i}$, there is $\bar{v}_{i} \in \mathcal{V}_{i}$ satisfying
\begin{equation}
\label{3-8}
  ( g_{i} \circ G_{i} ) ( \bar{x}, \bar{v}_{i} ) = \max\limits_{ v_{i} \in \mathcal{V}_{i} } ( g_{i} \circ G_{i} ) ( \bar{x}, v_{i} ) = \phi_{i} ( \bar{x} ).
\end{equation}
In addition, note that $\alpha \, \phi ( \bar{x} ) = 0$  due to $\alpha_{k} \, \phi ( x^{2k} ) \to 0$ as $k \to \infty$. Considering $\phi_{i} ( \bar{x} ) = \phi ( \bar{x} )$ for all $i \in I ( \bar{x} )$, we conclude from (\ref{3-8}) that
\begin{equation*}
  \mu_{i} \, ( g_{i} \circ G_{i} ) ( \bar{x}, \bar{v}_{i} ) = \dfrac{\alpha}{c} \, \bar{\mu}_{i} \, \phi_{i} ( \bar{x} ) =
  \dfrac{ \bar{\mu}_{i} }{c} \, [ \alpha \, \phi ( \bar{x} ) ] = 0,
\end{equation*}
i.e., $\mu_{i} \, ( g_{i} \circ G_{i} ) ( \bar{x}, \bar{v}_{i} ) = \mu_{i} \, \max\limits_{ v_{i} \in \mathcal{V}_{i} } ( g_{i} \circ G_{i} ) ( \bar{x}, v_{i} ) = 0$ for $i = 1, 2,\dots, n$. This yields the second relation of (\ref{3-1}).

\textbf{Case 2:} Next suppose that $\{ \alpha_{k} \}$ is unbounded. Similar as above, we have from $x^{*}_{2k} \in \alpha_{k} \, \widehat{\partial} \phi ( x^{2k} )$ that $( x^{*}_{2k}, - \alpha_{k} ) \in \widehat{N} ( ( x^{2k}, \phi ( x^{2k} ) ); \mathrm{gph}\,\phi )$ for each $k \in \mathbb{N}$. Hence
\begin{equation*}
  \Big( \dfrac{ x^{*}_{2k} }{ \alpha_{k} }, - 1 \Big) \in \widehat{N} ( ( x^{2k}, \phi ( x^{2k} ) ); \mathrm{gph}\,\phi ), \quad k \in \mathbb{N}.
\end{equation*}
Passing $k \to \infty$ and employing again (\ref{2-2}) give us some $( 0, - 1 ) \in N ( ( \bar{x}, \phi ( \bar{x} ) ); \mathrm{gph}\,\phi )$, which is equivalent to $0 \in \partial \phi ( \bar{x} )$. To proceed further as in the proof of the Case 1, there exists $\mu := ( \mu_{1}, \mu_{2},\dots, \mu_{n} ) \in \mathbb{R}^{n}_{+} \setminus \{ 0 \}$, with $\| \mu \| = 1 $, such that
\begin{equation*}
  0 \in \sum_{i=1}^{n} \mu_{i} \, \mathrm{cl}^{*}\mathrm{co} \medPl \medcup \medBl \medcup\limits_{ v_{i}^{*} \in \partial_{x} g_{i} ( G_{i} ( \bar{x}, v_{i} ) ) } \partial_{x} \langle v_{i}^{*}, G_{i} \rangle ( \bar{x}, v_{i} ) \,\mid\, v_{i} \in \mathcal{V}_{i} ( \bar{x} ) \medBr \medPr.
\end{equation*}
Moreover, noticing the unboundedness of $\{ \alpha_{k} \}$ and applying $\alpha_{k} \, \phi ( x^{2k} ) \to 0$ as $k \to \infty$, we may take $\bar{v}_{i} \in \mathcal{V}_{i}$ such that $\mu_{i} \, ( g_{i} \circ G_{i} ) ( \bar{x}, \bar{v}_{i} ) = \mu_{i} \, \phi_{i} ( \bar{x} ) =  \mu_{i} \, \phi ( \bar{x} ) = 0$ for each $i = 1, 2,\dots, n$. So, (\ref{3-1}) holds by choosing $y^{*} := 0 \in K^{+}$.
\end{proof}
\begin{rem}
\label{Rem3-1}
Theorem \ref{Thm3-2} collapses to
\begin{itemize}
  \item [(i)] (\citealt[Theorem~3.2]{Saadati(2022)}) for the problem $(\hyperlink{UP}{\mathrm{UP}})$ with $X = W = Z$, $\mathcal{V}_{i} = \mathcal{U}_{i}$, $i = 1, 2, \dots, n$, and identical maps $F$ and $G$,
  \item [(ii)] (\citealt[Proposition~3.9]{Fakhar(2018)}) with $X = W = Z$, $Y := \mathbb{R}^{p}$, $K := \mathbb{R}^{p}_{+}$, $\mathcal{V}_{i} = \mathcal{U}_{i}$, $i = 1, 2, \dots, m$, and identical maps $F$ and $G$, and
  \item [(iii)] (\citealt[Theorem~3.3]{Chuong(2016)}) with $X = W = Z := \mathbb{R}^{n}$, $Y := \mathbb{R}^{m}$, $K := \mathbb{R}^{m}_{+}$, $\mathcal{V}_{i} = \mathcal{U}_{i}$, $i = 1, 2, \dots, l$, and identical maps $F$ and $G$.
\end{itemize}
\end{rem}

The following corollary provides a Fritz-John optimality condition for weakly robust efficient solutions of the uncertain multiobjective optimization problem $(\hyperlink{SUP}{\mathrm{SUP}})$.
\begin{cor}
\label{Cor3-3}
Let $\bar{x}$ be a weakly robust efficient solution of the problem $(\hyperlink{SUP}{\mathrm{SUP}})$. Then there exist $y^{*} \in K^{+}$, $\mu := ( \mu_{1}, \mu_{2},\dots, \mu_{n} ) \in \mathbb{R}^{n}_{+}$, $\sigma := ( \sigma_{1}, \sigma_{2},\dots, \sigma_{m} ) \in \mathbb{R}^{m}_{+}$, with $ \| y^{*} \| + \| \mu \| + \| \sigma \| \ne 0$, and $\bar{v}_{i} \in \mathcal{V}_{i}$, $i = 1, 2,\dots, n$, such that
\begin{equation}
\label{3-9}
  \left\{
    \begin{aligned}
      & 0 \in \partial \langle y^{*}, f \rangle ( \bar{x} ) + \sum_{i=1}^{n} \mu_{i} \, \mathrm{cl}^{*}\mathrm{co} \medPl \medcup \medBl \partial_{x} g_{i} ( \bar{x}, v_{i} ) \,\mid\, v_{i} \in \mathcal{V}_{i} ( \bar{x} ) \medBr \medPr \\
      & \hspace{2.5cm} + \sum_{j=1}^{m} \sigma_{j} \, \mathrm{cl}^{*}\mathrm{co} \medPl \medcup \medBl \partial_{x} h_{j} ( \bar{x}, v_{n+j} ) \medcup \partial_{x} ( - h_{j} ) ( \bar{x}, v_{n+j} ) \,\mid\, v_{n+j} \in \mathcal{V}_{n+j} \medBr \medPr, \\
      & \mu_{i} \, \max\limits_{ v_{i} \in \mathcal{V}_{i} } g_{i} ( \bar{x}, v_{i} ) = \mu_{i} \, g_{i} ( \bar{x}, \bar{v}_{i} ) = 0, \quad i = 1, 2,\dots, n.
    \end{aligned}
  \right.
\end{equation}
\end{cor}
\begin{proof}
Note that the problem $(\hyperlink{SUP}{\mathrm{SUP}})$ is a particular case of the composite uncertain multiobjective optimization problem (\hyperlink{CUP}{CUP}), where $g$ is given as in (\ref{1-1}). First let us define the real-valued functions $\phi_{n+j}$, $j = 1, 2, \dots, m$, and $\phi$ on $X$ via
      \begin{equation*}
        \phi_{n+j} (x) := \max_{ v_{n+j} \in \mathcal{V}_{n+j} } | h_{j} \circ G_{n+j} | ( x, v_{n+j} ) \qquad \text{and} \qquad \phi (x) := \max_{ \substack{i = 1, 2,\dots, n \\ j = 1, 2, \dots, m} } \{ \phi_{i} (x), \phi_{n+j} (x) \}.
      \end{equation*}
Then proceeding similarly to the proof of Theorem \ref{Thm3-2} and invoking both the subdifferential chain rule
\begin{equation*}
\partial_{x} | h_{j} \circ G_{n+j} | ( x, v_{n+j}) \subset \hspace{-2mm} \medcup\limits_{-1 \le \tau \le 1} \hspace{-2mm} \partial_{x} ( \tau \, h_{j} \circ G_{n+j} ) ( x, v_{n+j} )
\end{equation*}
and the inclusion
\begin{equation*}
\partial_{x} ( \tau \, h_{j} \circ G_{n+j} ) ( x, v_{n+j} ) \subset | \tau |  ( \partial_{x} ( h_{j} \circ G_{n+j} ) ( x, v_{n+j} ) \medcup \partial_{x} -\!( h_{j} \circ G_{n+j} ) ( x, v_{n+j} ) ),
\end{equation*}
we find $y^{*} \in K^{+}$, $\gamma := ( \mu_{1}, \mu_{2},\dots, \mu_{n}, \sigma_{1}, \sigma_{2},\dots, \sigma_{m} ) \in \mathbb{R}^{n+m}_{+}$, with $ \| y^{*} \| + \| \gamma \| = 1$, $\bar{v}_{i} \in \mathcal{V}_{i}$, $i = 1, 2, \dots, n$, and $\bar{v}_{n+j} \in \mathcal{V}_{n+j}$, $j = 1, 2, \dots, m$, such that
\begin{equation*}
  \left\{
    \begin{aligned}
      & 0 \in \medcup\limits_{ w^{*} \in \partial \langle y^{*}, f \rangle ( F ( \bar{x} ) ) } \partial \langle w^{*}, F \rangle ( \bar{x} ) + \sum_{i=1}^{n} \mu_{i} \, \mathrm{cl}^{*}\mathrm{co} \medPl \medcup \medBl \medcup\limits_{ v_{i}^{*} \in \partial_{x} g_{i} ( G_{i} ( \bar{x}, v_{i} ) ) } \partial_{x} \langle v_{i}^{*}, G_{i} \rangle ( \bar{x}, v_{i} ) \\
      & \hspace{2mm} \,\mid\, v_{i} \in \mathcal{V}_{i} ( \bar{x} ) \medBr \medPr + \sum_{j=1}^{m} \sigma_{j} \, \mathrm{cl}^{*}\mathrm{co} \medPl \medcup \medBl \medcup\limits_{ v_{n+j}^{*} \in \partial_{x} h_{j} ( G_{n+j} ( \bar{x}, v_{n+j} ) ) \bigcup \partial_{x} - ( h_{j} ( G_{n+j} ( \bar{x}, v_{n+j} ) ) ) } \\
      & \hspace{2mm} \partial_{x} \langle v_{n+j}^{*}, G_{n+j} \rangle ( \bar{x}, v_{n+j} ) \,\mid\, v_{n+j} \in \mathcal{V}_{n+j} ( \bar{x} ) \medBr \medPr, \\
      & \mu_{i} \, \max\limits_{ v_{i} \in \mathcal{V}_{i} } g_{i} ( G_{i} ( \bar{x}, v_{i} ) ) = \mu_{i} \, g_{i} ( G_{i} ( \bar{x}, \bar{v}_{i} ) ) = 0, \quad i = 1, 2,\dots, n, \\
      & \sigma_{j} \, \max\limits_{ v_{n+j} \in \mathcal{V}_{n+j} } h_{j} ( G_{n+j} ( \bar{x}, v_{n+j} ) ) = \sigma_{j} \, h_{j} ( G_{n+j} ( \bar{x}, \bar{v}_{n+j} ) ) = 0, \quad j = 1, 2,\dots, m.
    \end{aligned}
  \right.
\end{equation*}
In this setting, we see that $X = W = Z$, $\mathcal{V}_{i} = \mathcal{U}_{i}$, $i = 1, 2, \dots, n$, $\mathcal{V}_{n+j} = \mathcal{U}_{n+j}$, $j = 1, 2, \dots, m$, and $F$ and $G$ are identical maps, thus the above relations reduces to the following ones
\begin{equation*}
  \left\{
    \begin{aligned}
      & 0 \in \partial \langle y^{*}, f \rangle ( \bar{x} ) + \sum_{i=1}^{n} \mu_{i} \, \mathrm{cl}^{*}\mathrm{co} \medPl \medcup \medBl \partial_{x} g_{i} ( \bar{x}, v_{i} ) \,\mid\, v_{i} \in \mathcal{V}_{i} ( \bar{x} ) \medBr \medPr \\
      & \hspace{2.5cm} + \sum_{j=1}^{m} \sigma_{j} \, \mathrm{cl}^{*}\mathrm{co} \medPl \medcup \medBl \partial_{x} h_{j} ( \bar{x}, v_{n+j} ) \medcup \partial_{x} ( - h_{j} ) ( \bar{x}, v_{n+j} ) \,\mid\, v_{n+j} \in \mathcal{V}_{n+j} \medBr \medPr, \\
      & \mu_{i} \, \max\limits_{ v_{i} \in \mathcal{V}_{i} } g_{i} ( \bar{x}, v_{i} ) = \mu_{i} \, g_{i} ( \bar{x}, \bar{v}_{i} ) = 0, \quad i = 1, 2,\dots, n,
    \end{aligned}
  \right.
\end{equation*}
due to $h_{j} ( \bar{x}, v_{n+j} ) = 0$ for all $v_{n+j} \in \mathcal{V}_{n+j}$, $j = 1, 2,\dots, m$. Clearly, $\| y^{*} \| + \| ( \mu_{1}, \mu_{2}, \dots, \mu_{n} ) \| + \| ( \sigma_{1}, \sigma_{2}, \dots, \sigma_{m} ) \| \ne 0$, and so the proof is complete.
\end{proof}

We now present an example which illustrates Theorem \ref{Thm3-2} for a composite uncertain multiobjective optimization problem.
\begin{exa}
\label{Exa3-2}
Let $X := \mathbb{R}^{2}$, $W := \mathbb{R}^{2}$, $Y := \mathbb{R}^{3}$, $Z := \mathbb{R}^{2}$, $\mathcal{V}_{i} = \mathcal{U}_{i} := [ - 1, 1 ]$, $i = 1, 2$, and $K := \mathbb{R}^{3}_{+}$. Consider the following composite uncertain optimization problem
\begin{equation*}
\hypertarget{CUPexa3.2}{}
 ( \mathrm{CUP} ) \qquad \min\nolimits_{K} \,\,\, \medBl ( f \circ F ) (x) \,\mid\, x := ( x_{1}, x_{2} ) \in X, \,\, ( g_{i} \circ G_{i} ) ( x, v_{i} ) \le 0, \,\, i = 1, 2 \medBr,
\end{equation*}
where $F : X \to W$, $F := ( F_{1}, F_{2} )$, is defined by $F_{1} ( x_{1}, x_{2} ) = \dfrac{1}{2} x_{1}$ and $F_{2} ( x_{1}, x_{2} ) = x_{2} - 1$, $f : W \to Y$, $f := ( f_{1}, f_{2}, f_{3} )$, is given by
\begin{equation*}
  \left\{
    \begin{aligned}
      f_{1} ( w_{1}, w_{2} ) & = -2 w_{1} + | w_{2} |, \\
      f_{2} ( w_{1}, w_{2} ) & = \dfrac{1}{ | w_{1} | + 1 } - 3 w_{2} + 2, \\
      f_{3} ( w_{1}, w_{2} ) & = \dfrac{1}{ \sqrt{ | w_{1} | + 1 } } - | w_{2} - 1 | - 1,
    \end{aligned}
  \right.
\end{equation*}
$G_{i} : X \times \mathcal{V}_{i} \to Z \times \mathcal{U}_{i}$, $i = 1, 2$, $G := ( G_{1}, G_{2} )$, are defined by $G_{1} ( x_{1}, x_{2}, v_{1} ) = ( x_{1} + 1, x_{2}, v_{1} )$ and $G_{2} ( x_{1}, x_{2}, v_{2} ) = ( x_{1}, 2 x_{2}, v_{2} )$, and $g_{i} : Z \times \mathcal{U}_{i} \to \mathbb{R}$, $i = 1, 2$, $g := ( g_{1}, g_{2} )$, are given by
\begin{equation*}
  \left\{
    \begin{aligned}
      g_{1} ( z_{1}, z_{2}, u_{1} ) & = u_{1}^{2} | z_{2} | + \max \big\{ z_{1}, 2 z_{1} \big\} - 3 | u_{1} |,\\
      g_{2} ( z_{1}, z_{2}, u_{2} ) & = -3 | z_{1} | + u_{2} z_{2} - 2,
    \end{aligned}
  \right.
\end{equation*}
where $v_{i} \in \mathcal{V}_{i}$ and $u_{i} \in \mathcal{U}_{i}$, $i = 1, 2$. It is easy to check that
\begin{align*}
  & \medBl ( x_{1}, x_{2} ) \in X \,\mid\, ( g_{1} \circ G_{1} ) ( x_{1}, x_{2}, v_{1} ) \le 0, \,\, \forall v_{1} \in \mathcal{V}_{1} \medBr \\
  = & \medBl ( x_{1}, x_{2} ) \in X \,\mid\, v_{1}^{2} | x_{2} | + \max \big\{ x_{1} + 1, 2 x_{1} + 2 \big\} - 3 | v_{1} | \le 0, \,\, \forall v_{1} \in \mathcal{V}_{1} \medBr \\
  = & \medBl ( x_{1}, x_{2} ) \in X \,\mid\, x_{1} \le - 1 \text{ and } | x_{2} | \le - x_{1} + 2 \medBr,
\end{align*}
and, since $x_{1} \le - 1$, it can be obtained that
\begin{align*}
    & \medBl ( x_{1}, x_{2} ) \in X \,\mid\, ( g_{2} \circ G_{2} ) ( x_{1}, x_{2}, v_{2} ) \le 0, \,\, \forall v_{2} \in \mathcal{V}_{2} \medBr \\
  = & \medBl ( x_{1}, x_{2} ) \in X \,\mid\, - 3 | x_{1} | + 2 v_{2} x_{2} - 2 \le 0, \,\, \forall v_{2} \in \mathcal{V}_{2} \medBr \\
  = & \medBl ( x_{1}, x_{2} ) \in X \,\mid\, x_{1} \le - 1 \text{ and } | x_{2} | \le - \dfrac{3}{2} x_{1} + 1 \medBr.
\end{align*}
Hence, the robust feasible set is
\begin{align*}
  C = & \medBl ( x_{1}, x_{2} ) \in X \,\mid\, - 2 \le x_{1} \le - 1 \text{ and } | x_{2} | \le - \dfrac{3}{2} x_{1} + 1 \medBr \medcup \\
      & \medBl ( x_{1}, x_{2} ) \in X \,\mid\, x_{1} \le - 2 \text{ and } | x_{2} | \le - x_{1} + 2 \medBr.
\end{align*}

Suppose that $\bar{x} := ( - 1, 1 ) \in C$ and $x := ( x_{1}, x_{2} ) \in C$. Considering $x_{1} \le - 1$, we have $( f_{1} \circ F_{1} ) (x) - ( f_{1} \circ F_{1} ) ( \bar{x} ) \ge 0$. Therefore
\begin{equation*}
  ( f \circ F ) (x) - ( f \circ F ) ( \bar{x} ) \notin - \mathrm{int}\hspace{.4mm}K, \quad \forall x \in C,
\end{equation*}
which means that $\bar{x}$ is a weakly robust efficient solution of the problem $(\hyperlink{CUPexa3.2}{\mathrm{CUP}})$. Observe also that
\begin{align*}
  \phi_{1} ( \bar{x} ) & = \max\limits_{ v_{1} \in \mathcal{V}_{1} } ( g_{1} \circ G_{1} ) ( \bar{x}, v_{1} ) = \max\limits_{ v_{1} \in \mathcal{V}_{1} } ( v_{1}^{2} - 3 | v_{1} | ) = 0, \\
  \phi_{2} ( \bar{x} ) & = \max\limits_{ v_{2} \in \mathcal{V}_{2} } ( g_{2} \circ G_{2} ) ( \bar{x}, v_{2} ) = \max\limits_{ v_{2} \in \mathcal{V}_{2} } ( 2 v_{2} - 5 ) = - 3.
\end{align*}
Hence $\phi ( \bar{x} ) = \max \big\{ \phi_{1} ( \bar{x} ), \phi_{2} ( \bar{x} ) \big\} = 0$, $\mathcal{V}_{1} ( \bar{x} ) = \{ 0 \}$, and $\mathcal{V}_{2} ( \bar{x} ) = \{ 1 \}$. Performing elementary calculations gives us $\partial F ( \bar{x} ) = \begin{pmatrix} \dfrac{1}{2} & 0 \\ 0 & 1 \end{pmatrix}$ and
\begin{equation*}
  \partial ( f_{1} \circ F_{1} ) ( \bar{x} ) = \{ - 1 \} \times [ - 1, 1 ], \, \partial ( f_{2} \circ F_{2} ) ( \bar{x} ) = \big[ - \dfrac{1}{2}, \dfrac{1}{2} \big] \times \{ - 3 \}, \, \partial ( f_{3} \circ F_{3} ) ( \bar{x} ) = \big[ - \dfrac{ \sqrt{2} }{2}, \dfrac{ \sqrt{2} }{2} \big] \times \{ - 1, 1 \}.
\end{equation*}
Further, we get $\partial_{x} G_{1} ( \bar{x}, v_{1} = 0 ) = \begin{pmatrix} 1 & 0 \\ 0 & 1 \end{pmatrix}$, $\partial_{x} G_{2} ( \bar{x}, v_{2} = 1 ) = \begin{pmatrix} 1 & 0 \\ 0 & 2 \end{pmatrix}$, and
\begin{equation*}
  \partial_{x} ( g_{1} \circ G_{1} ) ( \bar{x}, v_{1} = 0 ) = [ 1, 2 ] \times \{ 0 \}, \quad \partial_{x} ( g_{2} \circ G_{2} ) ( \bar{x}, v_{2} = 1 ) = \big\{ ( 3, 2 ), ( - 3, 2 ) \big\}.
\end{equation*}
So
\begin{equation}
\label{3-10}
  \left\{
    \begin{aligned}
  \mathrm{cl}^{*}\mathrm{co} \medPl \medcup\limits_{ v_{1}^{*} \in \partial_{x} g_{1} ( G_{1} ( \bar{x}, v_{1} ) ) } \partial_{x} \langle v_{1}^{*}, G_{1} \rangle ( \bar{x}, v_{1} = 0 ) \medPr & = [ 1, 2 ] \times \{ 0 \}, \\
  \mathrm{cl}^{*}\mathrm{co} \medPl \medcup\limits_{ v_{2}^{*} \in \partial_{x} g_{2} ( G_{2} ( \bar{x}, v_{2} ) ) } \partial_{x} \langle v_{2}^{*}, G_{2} \rangle ( \bar{x}, v_{2} = 1 ) \medPr & = [ - 3, 3 ] \times \{ 4 \}.
    \end{aligned}
  \right.
\end{equation}
On the other side, due to $I ( \bar{x} ) = \medBl i \in \{ 1, 2 \} \,\mid\, \phi_{i} ( \bar{x} ) = \phi ( \bar{x} ) \medBr = \{ 1 \}$, it is easy to see from (\ref{3-10}) that the $(\hyperlink{Def3-2}{\mathrm{CQ}})$, which will be presented later, is satisfied at $\bar{x}$.

Finally, we can find $y^{*} = \big( \dfrac{ \sqrt{2} }{3}, 0, \dfrac{ \sqrt{2} }{3} \big) \in K^{+} \setminus \{ 0 \}$ and $\mu = \big( \dfrac{1}{3}, 0 \big) \in \mathbb{R}^{2}_{+}$, with $\| y^{*} \| + \| \mu \| = 1$, satisfying
\begin{equation*}
  0 =
  \begin{pmatrix}
    \dfrac{\sqrt{2}}{3} & 0 & \dfrac{\sqrt{2}}{3}
  \end{pmatrix}
  \begin{pmatrix}
    - 1 &  0  & - \dfrac{1}{2} \\
      1 & - 3 & - 1
  \end{pmatrix}
  \begin{pmatrix}
    \dfrac{1}{2} & 0 \\
          0 & 1
  \end{pmatrix}
  +
  \begin{pmatrix}
    \dfrac{1}{3} & 0
  \end{pmatrix}
  \begin{pmatrix}
    \dfrac{3\sqrt{2}}{4} & 0 \\
               0 & 4
  \end{pmatrix},
\end{equation*}
and $\mu_{i} \, \max\limits_{ v_{i} \in \mathcal{V}_{i} } g_{i} ( G_{i} ( \bar{x}, v_{i} ) ) = 0$ for $i = 1, 2$.
\end{exa}

In order to establish sufficient optimality for (weakly) robust efficient solutions of the problem $(\hyperlink{CUP}{\mathrm{CUP}})$, we need to define a so-called robust \textit{Karush-Kuhn-Tucker} $(\hyperlink{Def3-1}{\mathrm{KKT}})$ condition for this problem.
\begin{defn}
\hypertarget{Def3-1}{}
A point $\bar{x} \in C$ is termed a robust (KKT) point if there exist $y^{*} \in K^{+} \setminus \{ 0 \}$, $\mu := ( \mu_{1}, \mu_{2},\dots, \mu_{n} ) \in \mathbb{R}^{n}_{+}$, and $\bar{v}_{i} \in \mathcal{V}_{i}$, $i=1, 2,\dots, n$, such that
\begin{equation*}
  \left\{
    \begin{aligned}
  & 0 \in \medcup\limits_{ w^{*} \in \partial \langle y^{*}, f \rangle ( F ( \bar{x} ) ) } \partial \langle w^{*}, F \rangle ( \bar{x} ) +
  \sum_{i=1}^{n} \mu_{i} \, \mathrm{cl}^{*}\mathrm{co} \medPl \medcup \medBl \medcup\limits_{ v_{i}^{*} \in \partial_{x} g_{i} ( G_{i} ( \bar{x}, v_{i} ) ) } \partial_{x} \langle v_{i}^{*}, G_{i} \rangle ( \bar{x}, v_{i} ) \,\mid\, v_{i} \in \mathcal{V}_{i} ( \bar{x} ) \medBr \medPr, \\
  & \mu_{i} \, \max\limits_{ v_{i} \in \mathcal{V}_{i} } ( g_{i} \circ G_{i} ) ( \bar{x}, v_{i} ) = \mu_{i} \, ( g_{i} \circ G_{i} ) ( \bar{x}, \bar{v}_{i} ) = 0, \quad i = 1, 2,\dots, n.
    \end{aligned}
  \right.
\end{equation*}
\end{defn}

It is worth to mention here that this condition $(\hyperlink{Def3-1}{\mathrm{KKT}})$ is in the sense of the Wolkowicz constraint qualification in the nonsmooth setting; see, e.g., \citealt{Wolkowicz(1980), Borwein(1986)} for more details. Moreover, it can be deduced from Theorem \ref{Thm3-2} that a weakly robust efficient solution of the problem $(\hyperlink{CUP}{\mathrm{CUP}})$ becomes a robust $(\hyperlink{Def3-1}{\mathrm{KKT}})$ point under the following \textit{constraint qualification} $(\hyperlink{Def3-2}{\mathrm{CQ}})$ condition in the sense of robustness.
\begin{defn}
\hypertarget{Def3-2}{}
{\normalfont (See (\citealt[Definition 2.3]{Saadati(2022)}))}
Let $\bar{x} \in C$. We say that the constraint qualification (CQ) is satisfied at $\bar{x}$ if
\begin{equation*}
  0 \notin \mathrm{cl}^{*}\mathrm{co} \medPl \medcup \medBl \medcup\limits_{ v_{i}^{*} \in \partial_{x} g_{i} ( G_{i} ( \bar{x}, v_{i} ) ) } \partial_{x} \langle v_{i}^{*}, G_{i} \rangle ( \bar{x}, v_{i} ) \,\mid\, v_{i} \in \mathcal{V}_{i} ( \bar{x} ) \medBr \medPr,\quad i \in I ( \bar{x} ),
\end{equation*}
where $I ( \bar{x} ) := \medBl i \in \{ 1, 2,\dots, n \} \,\mid\, \phi_{i} ( \bar{x} ) = \phi ( \bar{x} ) \medBr$.
\end{defn}

In general, a robust feasible point of the problem $(\hyperlink{CUP}{\mathrm{CUP}})$ at which the robust $(\hyperlink{Def3-1}{\mathrm{KKT}})$ condition holds may not be a (weakly) robust efficient solution. This motivates us to employ a similar concept of pseudo-quasi convexity in \citealt{Saadati(2022)} for the compositions $f \circ F$ and $g \circ G$.
\begin{defn}
\hypertarget{Def3-3}{}
\begin{itemize}
\item[(i)] We say that $( f \circ F, g \circ G )$ is \textit{type I pseudo convex} at $\bar{x} \in X$ if for any $x \in X$, $y^{*} \in K^{+}$, $w^{*} \in \partial \langle y^{*}, f \rangle ( F ( \bar{x} ) )$, $x^{*} \in \partial \langle w^{*}, F \rangle ( \bar{x} )$, $v^{*}_{i} \in \partial_{x} g_{i} ( G_{i} ( \bar{x}, v_{i} ) )$, and $x_{i}^{*} \in \partial_{x} \langle v^{*}_{i}, G_{i} \rangle ( \bar{x}, v_{i} )$, $v_{i} \in \mathcal{V}_{i} ( \bar{x} )$, $i = 1, 2, \dots, n$, there exists $\nu \in X$ such that
    \begin{align*}
      \langle y^{*}, f \circ F \rangle (x) < \langle y^{*}, f \circ F \rangle ( \bar{x} ) &\Longrightarrow \langle x^{*}, \nu \rangle  < 0, \\
      ( g_{i} \circ G_{i} ) ( x, v_{i} ) \le ( g_{i} \circ G_{i} ) ( \bar{x}, v_{i} ) &\Longrightarrow \langle x_{i}^{*}, \nu \rangle \le 0, \quad i = 1, 2,\dots, n.
    \end{align*}
\item[(ii)] We say that $( f \circ F, g \circ G )$ is \textit{type II pseudo convex} at $\bar{x} \in X$ if for any $x \in X \setminus \{ \bar{x} \}$, $y^{*} \in K^{+} \setminus \{ 0 \}$, $w^{*} \in \partial \langle y^{*}, f \rangle ( F ( \bar{x} ) )$, $x^{*} \in \partial \langle w^{*}, F \rangle ( \bar{x} )$, $v^{*}_{i} \in \partial_{x} g_{i} ( G_{i} ( \bar{x}, v_{i} ) )$, and $x_{i}^{*} \in \partial_{x} \langle v^{*}_{i}, G_{i} \rangle ( \bar{x}, v_{i} )$, $v_{i} \in \mathcal{V}_{i} ( \bar{x} )$, $i = 1, 2,\dots, n$, there exists $\nu \in X$ such that
    \begin{align*}
      \langle y^{*}, f \circ F \rangle (x) \le \langle y^{*}, f \circ F \rangle ( \bar{x} ) &\Longrightarrow \langle x^{*}, \nu \rangle  < 0, \\
      ( g_{i} \circ G_{i} ) ( x, v_{i} ) \le ( g_{i} \circ G_{i} ) ( \bar{x}, v_{i} ) &\Longrightarrow \langle x_{i}^{*}, \nu \rangle \le 0, \quad i = 1, 2,\dots, n.
    \end{align*}
\end{itemize}
\end{defn}
\begin{rem}
\label{Rem3-3}
In Definition \hyperlink{Def3-3}{3.3},
\begin{itemize}
  \item[(i)] if $X = W = Z$, $\mathcal{V}_{i} = \mathcal{U}_{i}$, $i = 1, 2, \dots, n$, and $F$ and $G$ are identical maps, then this definition reduces to the corresponding one in (\citealt[Definition~2.2]{Saadati(2022)}). Moreover, as shown in (\citealt[Example~2.2]{Saadati(2022)}), the class of \textit{type I pseudo convex} functions is properly larger than the class of \textit{generalized convex} functions.
  \item[(ii)] if $X = W = Z$, $Y := \mathbb{R}^{p}$, $K := \mathbb{R}^{p}_{+}$, $\mathcal{V}_{i} = \mathcal{U}_{i}$, $i = 1, 2, \dots, m$, and $F$ and $G$ are identical maps, then this definition collapses to (\citealt[Definition~3.2]{Fakhar(2018)}). Furthermore, as demonstrated in (\citealt[Example~3.4]{Fakhar(2018)}), the class of \textit{pseudo-quasi generalized convex} functions is properly wider than the class of generalized convex functions.
  \item[(iii)] if $X = W = Z := \mathbb{R}^{n}$, $Y := \mathbb{R}^{m}$, $K := \mathbb{R}^{m}_{+}$, $\mathcal{V}_{i} = \mathcal{U}_{i}$, $i = 1, 2, \dots, l$, and $F$ and $G$ are identical maps, then this definition reduces to (\citealt[Definition~3.9]{Chuong(2016)}). Furthermore, as illustrated in (\citealt[Example~3.10]{Chuong(2016)}), the class of generalized convex functions contains some \textit{nonconvex} functions.
\end{itemize}
\end{rem}

\begin{rem}
\label{Rem3-4}
It follows from Definition \hyperlink{Def3-3}{3.3} that if $( f \circ F, g \circ G )$ is type II pseudo convex at $\bar{x} \in X$, then $( f \circ F, g \circ G )$ is type I pseudo convex at $\bar{x} \in X$, but converse is not true (see (\citealt[Example~2.2]{Saadati(2022)}) considering special case of Remark \ref{Rem3-3}(i)).
\end{rem}

The forthcoming proposition shows that the class of (resp., type II) type I pseudo convex composite functions includes (resp., \textit{strictly}) \textit{convex} composite functions.
\begin{prop}
\label{Prop3-4}
Let $\bar{x} \in X$, and let $f$, $F$, $g$, and $G$ be such that $\langle y^{*}, f \rangle$ is convex for every $y^{*} \in K^{+}$, $\langle w^{*}, F \rangle$ is convex for every $w^{*} \in \partial \langle y^{*}, f \rangle ( F ( \bar{x} ) )$, $g_{i}$ is a convex function, and $\langle v^{*}_{i}, G_{i} \rangle$ is convex for every $v^{*}_{i} \in \partial_{x} g_{i} ( G_{i} ( \bar{x}, v_{i} ) )$, $v_{i} \in \mathcal{V}_{i} ( \bar{x} )$, $i = 1, 2, \dots, n$. Then
\begin{itemize}
  \item[{\normalfont (i)}] $( f \circ F, g \circ G )$ is type I pseudo convex at $\bar{x}$.
  \item[{\normalfont (ii)}] If $\langle y^{*}, f \rangle$ is strictly convex for every $y^{*} \in K^{+} \setminus \{ 0 \}$ and $F$ is injective, then $( f \circ F, g \circ G )$ is type II pseudo convex at $\bar{x}$.
\end{itemize}
\end{prop}
\begin{proof}
(i) Suppose that $x \in X$, $y^{*} \in K^{+}$, $w^{*} \in \partial \langle y^{*}, f \rangle ( F ( \bar{x} ) )$, $x^{*} \in \partial \langle w^{*}, F \rangle ( \bar{x} )$, $v^{*}_{i} \in \partial_{x} g_{i} ( G_{i} ( \bar{x}, v_{i} ) )$, and $x_{i}^{*} \in \partial_{x} \langle v^{*}_{i}, G_{i} \rangle ( \bar{x}, v_{i} )$, $v_{i} \in \mathcal{V}_{i} ( \bar{x} )$, $i = 1, 2, \dots, n$, and that
\begin{align}
\label{3-11}
  \langle y^{*}, f \circ F \rangle (x) & < \langle y^{*}, f \circ F \rangle ( \bar{x} ), \\
\label{3-12}
  ( g_{i} \circ G_{i} ) ( x, v_{i} ) & \le ( g_{i} \circ G_{i} ) ( \bar{x}, v_{i} ), \quad i = 1, 2,\dots, n.
\end{align}
Take $\nu := x - \bar{x}$. Since $\langle w^{*}, F \rangle$ is a convex function, we have
\begin{equation*}
  \langle x^{*}, \nu \rangle = \langle x^{*}, x - \bar{x} \rangle \le \langle w^{*}, F \rangle (x) - \langle w^{*}, F \rangle ( \bar{x} ) = \langle w^{*}, F (x) - F ( \bar{x} ) \rangle.
\end{equation*}
On the other hand, since $\langle y^{*}, f \rangle$ is convex, we obtain
\begin{equation*}
  \langle w^{*}, F (x) - F ( \bar{x} ) \rangle \le \langle y^{*}, f \rangle ( F (x) ) - \langle y^{*}, f \rangle ( F ( \bar{x} ) ) = \langle y^{*}, f \circ F \rangle (x) - \langle y^{*}, f \circ F \rangle ( \bar{x} ).
\end{equation*}
Combining now the both latter inequalities with (\ref{3-11}), one gets
\begin{equation*}
  \langle x^{*}, \nu \rangle \le \langle y^{*}, f \circ F \rangle (x) - \langle y^{*}, f \circ F \rangle ( \bar{x} ) < 0.
\end{equation*}
Similarly, using the convexity of functions $g_{i}$ and $\langle v^{*}_{i}, G_{i} \rangle$ and applying (\ref{3-12}) give us
\begin{equation*}
  \langle x^{*}_{i}, \nu \rangle \le ( g_{i} \circ G_{i} ) ( x, v_{i} ) - ( g_{i} \circ G_{i} ) ( \bar{x}, v_{i} ) \le 0, \quad i = 1, 2,\dots, n.
\end{equation*}
Consequently, $( f \circ F, g \circ G )$ is type I pseudo convex at $\bar{x}$.

(ii) Suppose that $x \in X \setminus \{ \bar{x} \}$, $y^{*} \in K^{+} \setminus \{ 0 \}$, $w^{*} \in \partial \langle y^{*}, f \rangle ( F ( \bar{x} ) )$, $x^{*} \in \partial \langle w^{*}, F \rangle ( \bar{x} )$, $v^{*}_{i} \in \partial_{x} g_{i} ( G_{i} ( \bar{x}, v_{i} ) )$, and $x_{i}^{*} \in \partial_{x} \langle v^{*}_{i}, G_{i} \rangle ( \bar{x}, v_{i} )$, $v_{i} \in \mathcal{V}_{i} ( \bar{x} )$, $i = 1, 2, \dots, n$, and that
\begin{align}
\label{3-13}
  \langle y^{*}, f \circ F \rangle (x) & \le \langle y^{*}, f \circ F \rangle ( \bar{x} ), \\
  ( g_{i} \circ G_{i} ) ( x, v_{i} ) & \le ( g_{i} \circ G_{i} ) ( \bar{x}, v_{i} ), \quad i = 1, 2,\dots, n. \nonumber
\end{align}
Choosing $\nu := x - \bar{x}$ and proceeding as in the part (i), it holds that
\begin{equation}
\label{3-14}
  \langle x^{*}, \nu \rangle \le \langle w^{*}, F (x) - F ( \bar{x} ) \rangle
\end{equation}
and
\begin{equation*}
  \langle x^{*}_{i}, \nu \rangle \le ( g_{i} \circ G_{i} ) ( x, v_{i} ) - ( g_{i} \circ G_{i} ) ( \bar{x}, v_{i} ) \le 0, \quad i = 1, 2,\dots, n.
\end{equation*}
Note here that $F (x) \ne F ( \bar{x} )$ due to the injectivity of $F$. Since $\langle y^{*}, f \rangle$ is strictly convex, it follows from (\citealt[Proposition~6.1.3]{Hiriart-Urruty(1993)B}) that
\begin{equation*}
  \langle w^{*}, F (x) - F ( \bar{x} ) \rangle < \langle y^{*}, f \rangle ( F (x) ) - \langle y^{*}, f \rangle ( F ( \bar{x} ) ) = \langle y^{*}, f \circ F \rangle (x) - \langle y^{*}, f \circ F \rangle ( \bar{x} ).
\end{equation*}
This combined with (\ref{3-13}) and (\ref{3-14}) yields that
\begin{equation*}
  \langle x^{*}, \nu \rangle < \langle y^{*}, f \circ F \rangle (x) - \langle y^{*}, f \circ F \rangle ( \bar{x} ) \le 0.
\end{equation*}
So, $( f \circ F, g \circ G )$ is type II pseudo convex at $\bar{x}$.
\end{proof}

Now we are ready to derive a robust $(\hyperlink{Def3-1}{\mathrm{KKT}})$ sufficient condition for (weakly) robust efficient solutions of the problem $(\hyperlink{CUP}{\mathrm{CUP}})$.
\begin{thm}
\label{Thm3-5}
Assume that $\bar{x} \in C$ satisfies the robust $(\hyperlink{Def3-1}{\mathrm{KKT}})$ condition.
\begin{itemize}
  \item[{\normalfont (i)}] If $( f \circ F, g \circ G )$ is type I pseudo convex at $\bar{x}$, then $\bar{x} \in \mathcal{S}^{w}(\hyperlink{CRP}{\mathrm{CRP}})$.
  \item[{\normalfont (ii)}] If $( f \circ F, g \circ G )$ is type II pseudo convex at $\bar{x}$, then $\bar{x} \in \mathcal{S}(\hyperlink{CRP}{\mathrm{CRP}})$.
\end{itemize}
\end{thm}
\begin{proof}
As $\bar{x} \in C$ is a robust $(\hyperlink{Def3-1}{\mathrm{KKT}})$ point of the problem $(\hyperlink{CUP}{\mathrm{CUP}})$, we can find $y^{*} \in K^{+} \setminus \{ 0 \}$, $w^{*} \in \partial \langle y^{*}, f \rangle ( F ( \bar{x} ) )$, $x^{*} \in \partial \langle w^{*}, F \rangle ( \bar{x} )$, $\mu_{i} \ge 0$, and $u_{i}^{*} \in \mathrm{cl}^{*}\mathrm{co} \medPl \medcup \medBl \medcup\limits_{ v_{i}^{*} \in \partial_{x} g_{i} ( G_{i} ( \bar{x}, v_{i} ) ) } \partial_{x} \langle v_{i}^{*}, G_{i} \rangle ( \bar{x}, v_{i} ) \,\mid\, v_{i} \in \mathcal{V}_{i} ( \bar{x} ) \medBr \medPr$, $i = 1, 2,\dots, n$, such that
\begin{align}
\label{3-15}
  & 0 = x^{*} + \sum\limits_{i=1}^{n} \mu_{i} \, u_{i}^{*}, \\
\label{3-16}
  & \mu_{i} \, \max\limits_{ v_{i} \in \mathcal{V}_{i} } g_{i} ( G_{i} ( \bar{x}, v_{i} ) ) = 0, \quad i = 1, 2,\dots, n.
\end{align}
First prove (i). On the contrary, suppose that $\bar{x} \notin \mathcal{S}^{w}(\hyperlink{CRP}{\mathrm{CRP}})$. This means by definition that there exists $\hat{x} \in C$ such that $( f \circ F ) ( \hat{x} ) - ( f \circ F ) ( \bar{x} ) \in - \mathrm{int}\hspace{.4mm}K$. It follows from (\citealt[Lemma~3.21]{Jahn(2004)B}) that $\langle y^{*}, ( f \circ F ) ( \hat{x} ) - ( f \circ F ) ( \bar{x} ) \rangle < 0$. Noting that $( f \circ F, g \circ G )$ is the type I pseudo convex at $\bar{x}$, we conclude from the latter inequality that there is $\nu \in X$ such that
\begin{equation}
\label{3-17}
  \langle x^{*}, \nu \rangle < 0.
\end{equation}
In addition, taking (\ref{3-15}) into account gives us
\begin{equation}
\label{3-18}
  0 = \langle x^{*}, \nu \rangle + \sum_{i=1}^{n} \mu_{i} \, \langle u_{i}^{*}, \nu \rangle
\end{equation}
for $\nu$ above. So, the relationships in (\ref{3-17}) and (\ref{3-18}) lead to
\begin{equation*}
  \sum_{i=1}^{n} \mu_{i} \, \langle u_{i}^{*}, \nu \rangle > 0.
\end{equation*}
To proceed, first pick $i_{0} \in \{ 1, 2,\dots, n \}$ satisfying $\mu_{ i_{0} } \, \langle u_{ i_{0} }^{*}, \nu \rangle > 0$. On the one side, due to
\begin{equation*}
u_{ i_{0} }^{*} \in \mathrm{cl}^{*}\mathrm{co} \medPl \medcup \medBl \medcup\limits_{ v_{ i_{0} }^{*} \in \partial_{x} g_{ i_{0} } ( G_{ i_{0} } ( \bar{x}, v_{ i_{0} } ) ) } \partial_{x} \langle v_{ i_{0} }^{*}, G_{ i_{0} } \rangle ( \bar{x}, v_{ i_{0} } ) \,\mid\, v_{ i_{0} } \in \mathcal{V}_{ i_{0} } ( \bar{x} ) \medBr \medPr,
\end{equation*}
we have a sequence $\{ u_{ i_{0} k }^{*} \} \subset \mathrm{co} \medPl \medcup \medBl \medcup\limits_{ v_{ i_{0} }^{*} \in \partial_{x} g_{ i_{0} } ( G_{ i_{0} } ( \bar{x}, v_{ i_{0} } ) ) } \partial_{x} \langle v_{ i_{0} }^{*}, G_{ i_{0} } \rangle ( \bar{x}, v_{ i_{0} } ) \,\mid\, v_{ i_{0} } \in \mathcal{V}_{ i_{0} } ( \bar{x} ) \medBr \medPr$ such that $u_{ i_{0} k }^{*} \overset{\scriptscriptstyle w^{*}} \to u_{ i_{0} }^{*}$. Since $\mu_{ i_{0} } > 0$, there exists $k_{0} \in \mathbb{N}$ such that
\begin{equation}
\label{3-19}
  \langle u_{ i_{0} k_{0} }^{*}, \nu \rangle > 0.
\end{equation}
On the other side, due to $u_{ i_{0} k_{0} }^{*} \in \mathrm{co} \medPl \medcup \medBl \medcup\limits_{ v_{ i_{0} }^{*} \in \partial_{x} g_{ i_{0} } ( G_{ i_{0} } ( \bar{x}, v_{ i_{0} } ) ) } \partial_{x} \langle v_{ i_{0} }^{*}, G_{ i_{0} } \rangle ( \bar{x}, v_{ i_{0} } ) \,\mid\, v_{ i_{0} } \in \mathcal{V}_{ i_{0} } ( \bar{x} ) \medBr \medPr$, we get $u_{p}^{*} \in \medcup \medBl \medcup\limits_{ v_{ i_{0} }^{*} \in \partial_{x} g_{ i_{0} } ( G_{ i_{0} } ( \bar{x}, v_{ i_{0} } ) ) } \partial_{x} \langle v_{ i_{0} }^{*}, G_{ i_{0} } \rangle ( \bar{x}, v_{ i_{0} } ) \,\mid\, v_{ i_{0} } \in \mathcal{V}_{ i_{0} } ( \bar{x} ) \medBr$ and $\mu_{p} \ge 0$ with $\mathlarger{\sum\limits}_{p=1}^{s} \, \mu_{p} = 1$, $p = 1, 2,\dots, s$, $s \in \mathbb{N}$, satisfying $u_{ i_{0} k_{0} }^{*} = \mathlarger{\sum\limits}_{p=1}^{s} \, \mu_{p} \, u_{p}^{*}$. This combined with (\ref{3-19}) entails that $\mathlarger{\sum\limits}_{p=1}^{s} \, \mu_{p} \, \langle u_{p}^{*}, \nu \rangle > 0$. Then we may choose $p_{0} \in \{ 1, 2,\dots, s \}$ so that
\begin{equation}
\label{3-20}
  \langle u_{ p_{0} }^{*}, \nu \rangle > 0,
\end{equation}
and take $\bar{v}_{ i_{0} } \in \mathcal{V}_{ i_{0} } ( \bar{x} )$ and $\bar{v}_{ i_{0} }^{*} \in \partial_{x} g_{ i_{0} } ( G_{ i_{0} } ( \bar{x}, \bar{v}_{ i_{0} } ) )$ such that $u_{ p_{0} }^{*} \in \partial_{x} \langle \bar{v}_{ i_{0} }^{*}, G_{ i_{0} } \rangle ( \bar{x}, \bar{v}_{ i_{0} } )$ due to $u_{ p_{0} }^{*} \in \medcup \medBl \medcup\limits_{ v_{ i_{0} }^{*} \in \partial_{x} g_{ i_{0} } ( G_{ i_{0} } ( \bar{x}, v_{ i_{0} } ) ) } \partial_{x} \langle v_{ i_{0} }^{*}, G_{ i_{0} } \rangle ( \bar{x}, v_{ i_{0} } ) \,\mid\, v_{ i_{0} } \in \mathcal{V}_{ i_{0} } ( \bar{x} ) \medBr$. Employing now the type I pseudo convexity of $( f \circ F, g \circ G )$ at $\bar{x}$ and applying (\ref{3-20}), one has
\begin{equation}
\label{3-21}
  ( g_{ i_{0} } \circ G_{ i_{0} } ) ( \hat{x}, \bar{v}_{ i_{0} } ) > ( g_{ i_{0} } \circ G_{ i_{0} } ) ( \bar{x}, \bar{v}_{ i_{0} } ).
\end{equation}
Since $\bar{v}_{ i_{0} } \in \mathcal{V}_{ i_{0} } ( \bar{x} )$, it implies that $g_{ i_{0} } ( G_{ i_{0} } ( \bar{x}, \bar{v}_{ i_{0} } ) ) = \max\limits_{ v_{ i_{0} } \in \mathcal{V}_{ i_{0} } } g_{ i_{0} } ( G_{ i_{0} } ( \bar{x}, v_{ i_{0} } ) )$ which gives by (\ref{3-16}) that $\mu_{ i_{0} } \, g_{ i_{0} } ( G_{ i_{0} } ( \bar{x}, \bar{v}_{ i_{0} } ) ) = 0$. Using the last equality together with (\ref{3-21}) allows us $\mu_{ i_{0} } \, ( g_{ i_{0} } \circ G_{ i_{0} } ) ( \hat{x}, \bar{v}_{ i_{0} } ) > 0$, and so $( g_{ i_{0} } \circ G_{ i_{0} } ) ( \hat{x}, \bar{v}_{ i_{0} } ) > 0$, which contradicts with the fact that $\hat{x} \in C$ and ends the proof of (i).

The verification of assertion (ii) is similar to the part (i). Let $\bar{x} \notin \mathcal{S}(\hyperlink{CRP}{\mathrm{CRP}})$. Then there exists $\hat{x} \in C$ satisfying $( f \circ F ) ( \hat{x} ) - ( f \circ F ) ( \bar{x} ) \in - K \setminus \{ 0 \}$. This holds that $\hat{x} \ne \bar{x}$ and $\langle y^{*}, ( f \circ F ) ( \hat{x} ) - ( f \circ F ) ( \bar{x} ) \rangle \le 0$. Finally, involving the type II pseudo convexity of $( f \circ F, g \circ G )$  at $\bar{x}$, we arrive at the result.
\end{proof}
\begin{rem}
\label{Rem3-5}
Theorem \ref{Thm3-5} reduces to
\begin{itemize}
  \item [(i)] (\citealt[Theorem~3.4]{Saadati(2022)}) under pseudo-quasi convex assumptions,
  \item [(ii)] (\citealt[Theorem~3.10]{Fakhar(2018)}) where the involved functions are pseudo-quasi generalized convexity, and
  \item [(iii)] (\citealt[Theorem~3.11]{Chuong(2016)}) with generalized convex functions,
\end{itemize}
when the last two deal with the finite-dimensional frameworks.
\end{rem}

The following corollary of Theorem \ref{Thm3-5} concerns a convex problem of uncertain multiobjective optimization to reobtain the robust $(\hyperlink{Def3-1}{\mathrm{KKT}})$ sufficient optimality for (weakly) robust efficient solutions.
\begin{cor}
\label{Cor3-6}
Let $f$ and $g_{i}$, $i = 1, 2,\dots, n$, be convex functions and $h_{j}$, $j = 1, 2,\dots, m$, be affine functions. Suppose that $\bar{x}$ is a robust $(\hyperlink{Def3-1}{\mathrm{KKT}})$ point of the problem $(\hyperlink{SUP}{\mathrm{SUP}})$, i.e., $(\mathrm{\ref{3-9}})$ holds with $y^{*} \ne 0$. Then $\bar{x}$ is a weakly robust efficient solution of the problem $(\hyperlink{SUP}{\mathrm{SUP}})$. If $f$ is a strictly convex function, then $\bar{x}$ is a robust efficient solution of such problem.
\end{cor}
\begin{proof}
Observe, as in the proof of Corollary \ref{Cor3-3}, that the problem $(\hyperlink{SUP}{\mathrm{SUP}})$ is a special case of the composite uncertain multiobjective optimization problem $(\hyperlink{CUP}{\mathrm{CUP}})$, where $g$ is given as in (\ref{1-1}). Invoking the convexity of $f$ and $g_{i}$ and the affineness of $h_{j}$, it follows from Proposition \ref{Prop3-4}(i) that $( f \circ F, g \circ G )$ is type II pseudo convex at $\bar{x}$. If suppose that $f$ is a strict convex function, then $( f \circ F, g \circ G )$ is type II pseudo convex at $\bar{x}$ as shown by Proposition \ref{Prop3-4}(ii). So, we directly arrive at the desired conclusions by applying Theorem \ref{Thm3-5}.
\end{proof}
\section{Robust duality}
\label{Sec4-Dual}
In this section, we address a \textit{Mond-Weir-type dual} problem for the composite robust multiobjective optimization problem $(\hyperlink{CRP}{\mathrm{CRP}})$, and study the weak, strong, and converse duality relations between the corresponding problems under the assumptions of pseudo-quasi convexity.

Let $z \in X$, $y^{*} \in K^{+} \setminus \{ 0 \}$, and $\mu \in \mathbb{R}^{n}_{+}$. In connection with the problem $(\hyperlink{CRP}{\mathrm{CRP}})$, we introduce a \textit{dual robust} problem of the form
\begin{equation*}
\hypertarget{CRD}{}
  ( \mathrm{CRD} ) \qquad \max\nolimits_{K} \,\,\, \medBl \bar{f} ( z, y^{*}, \mu ) := ( f \circ F ) (z) \,\mid\, ( z, y^{*}, \mu ) \in C_{D} \medBr,
\end{equation*}
where $C_{D}$ is the feasible set defined by
\begin{align*}
  C_{D} := & \Big\{ ( z, y^{*}, \mu ) \in \,\, X \times ( K^{+} \setminus \{ 0 \} ) \times \mathbb{R}^{n}_{+} \,\mid\, 0 \in \medcup\limits_{ w^{*} \in \partial \langle y^{*}, f \rangle ( F ( \bar{x} ) ) } \partial \langle w^{*}, F \rangle ( \bar{x} ) \\
  & + \sum_{i=1}^{n} \mu_{i} \, \mathrm{cl}^{*}\mathrm{co} \medPl \medcup \medBl \medcup\limits_{ v_{i}^{*} \in \partial_{x} g_{i} ( G_{i} ( \bar{x}, v_{i} ) ) } \partial_{x} \langle v_{i}^{*}, G_{i} \rangle ( \bar{x}, v_{i} ) \,\mid\, v_{i} \in \mathcal{V}_{i} ( \bar{x} ) \medBr \medPr, \,\, \\
  & \mu_{i} \, g_{i} ( G_{i} ( z, v_{i} ) ) \ge 0, \,\, i = 1, 2,\dots, n \Big\}.
\end{align*}

It should be noticed here that the notions of robust efficient solutions (resp., weakly robust efficient solutions) of the dual problem $(\hyperlink{CRD}{\mathrm{CRD}})$ are understood as in Definition \hyperlink{Def1-1}{1.1} by replacing $-K$ (resp., $-\mathrm{int}\hspace{.4mm}K$) by $K$ (resp., $\mathrm{int}\hspace{.4mm}K$). Besides, denote by $\mathcal{S}(\hyperlink{CRD}{\mathrm{CRD}})$ (resp., $\mathcal{S}^{w}(\hyperlink{CRD}{\mathrm{CRD}})$) the set of robust efficient solutions (resp., weakly robust efficient solutions) of the problem $(\hyperlink{CRD}{\mathrm{CRD}})$. And make for convenience the standard notations
\begin{align*}
  u \prec v \Leftrightarrow u-v \in -\mathrm{int}\hspace{.4mm}K, \quad & u \nprec v \text{ is the negation of } u \prec v, \\
  u \preceq v \Leftrightarrow u-v \in -K \setminus\{ 0 \}, \quad & u \npreceq v \text{ is the negation of } u \preceq v.
\end{align*}

Weak duality theorem which holds between the primal problem $(\hyperlink{CRP}{\mathrm{CRP}})$ and the dual problem $(\hyperlink{CRD}{\mathrm{CRD}})$ is represented as follows.
\begin{thm}
\label{Thm4-1}
\normalfont{\textbf{\textsc{(Weak Duality)}}}
Let $x \in C$, and let $( z, y^{*}, \mu ) \in C_{D}$.
\begin{itemize}
  \item[{\normalfont (i)}] If $( f \circ F, g \circ G )$ is type I pseudo convex at $z$, then $( f \circ F ) (x) \nprec \bar{f} ( z, y^{*}, \mu )$.
  \item[{\normalfont (ii)}] If $( f \circ F, g \circ G )$ is type II pseudo convex at $z$, then $( f \circ F ) (x) \npreceq \bar{f} ( z, y^{*}, \mu )$.
  \end{itemize}
\end{thm}
\begin{proof}
Since $( z, y^{*}, \mu ) \in C_{D}$, there exist $y^{*} \in K^{+} \setminus \{ 0 \}$, $w^{*} \in \partial \langle y^{*}, f \rangle ( F ( \bar{x} ) )$, $x^{*} \in \partial \langle w^{*}, F \rangle ( \bar{x} )$, $\mu_{i} \ge 0$, and $u_{i}^{*} \in \mathrm{cl}^{*}\mathrm{co} \medPl \medcup \medBl \medcup\limits_{ v_{i}^{*} \in \partial_{x} g_{i} ( G_{i} ( \bar{x}, v_{i} ) ) } \partial_{x} \langle v_{i}^{*}, G_{i} \rangle ( \bar{x}, v_{i} ) \,\mid\, v_{i} \in \mathcal{V}_{i} ( \bar{x} ) \medBr \medPr$, $i = 1, 2,\dots, n$, satisfying
\begin{align}
\label{4-1}
  & 0 = x^{*} + \sum\limits_{i=1}^{n} \mu_{i} \, u_{i}^{*}, \\
  & \mu_{i} \, \max\limits_{ v_{i} \in \mathcal{V}_{i} } g_{i} ( G_{i} ( \bar{x}, v_{i} ) ) = 0, \quad i = 1, 2,\dots, n. \nonumber
\end{align}

To prove (i) by contradiction, suppose that $( f \circ F ) (x) \prec \bar{f} ( z, y^{*}, \mu )$. Then $\langle y^{*}, ( f \circ F ) (x) - \bar{f} ( z, y^{*}, \mu ) \rangle < 0$ due to $y^{*} \ne 0$. This is nothing else but $\langle y^{*}, ( f \circ F ) (x) - ( f \circ F ) (z) \rangle < 0$. Using the type I pseudo convexity of $( f \circ F, g \circ G )$ at $z$, we infer from the latter inequality that there is $\nu \in X$ such that
\begin{equation*}
  \langle x^{*}, \nu \rangle < 0.
\end{equation*}
On the other side, applying (\ref{4-1}) allows us
\begin{equation*}
  0 = \langle x^{*}, \nu \rangle + \sum_{i=1}^{n} \mu_{i} \, \langle u_{i}^{*}, \nu \rangle
\end{equation*}
for $\nu$ above. Combining these two relations, we arrive at
\begin{equation*}
  \sum_{i=1}^{n} \mu_{i} \, \langle x^{*}_{i}, \nu \rangle > 0.
\end{equation*}
Now assume that there exists $i_{0} \in \{ 1, 2,\dots, n \}$ such that $\mu_{ i_{0} } \, \langle x^{*}_{ i_{0} }, \nu \rangle > 0$. Proceeding as in the proof of Theorem \ref{Thm3-5}(i) and changing from $\hat{x} - \bar{x}$ to $x - z$ yield $( g_{ i_{0} } \circ G_{ i_{0} } ) ( x, \bar{v}_{ i_{0} } ) > 0$, which contradicts with $x \in C$.

Next to justify (ii), we follow the proof of the part (i) by definition of the type II pseudo convexity of $( f \circ F, g \circ G )$ at $z$ and observe that the condition $( f \circ F ) (x) \preceq \bar{f} ( z, y^{*}, \mu )$ leads us to a contradiction.
\end{proof}

The two forthcoming theorems declare strong duality relationships between the primal problem $(\hyperlink{CRP}{\mathrm{CRP}})$ and the dual problem $(\hyperlink{CRD}{\mathrm{CRD}})$.
\begin{thm}
\label{Thm4-2}
\normalfont{\textbf{\textsc{(Strong Duality)}}}
Let $\bar{x} \in \mathcal{S}^{w}(\hyperlink{CRP}{\mathrm{CRP}})$ be such that the $(\hyperlink{Def3-2}{\mathrm{CQ}})$ is satisfied at this point. Then there exists $( \bar{y}^{*}, \bar{\mu} ) \in K^{+} \setminus \{ 0 \} \times \mathbb{R}^{n}_{+}$ such that $( \bar{x}, \bar{y}^{*}, \bar{\mu} ) \in C_{D}$. Furthermore,
\begin{itemize}
  \item[{\normalfont (i)}] If $( f \circ F, g \circ G )$ is type I pseudo convex at $z$ for all $z \in X$, then $( \bar{x}, \bar{y}^{*}, \bar{\mu} ) \in \mathcal{S}^{w}(\hyperlink{CRD}{\mathrm{CRD}})$.
  \item[{\normalfont (ii)}] If $( f \circ F, g \circ G )$ is type II pseudo convex at $z$ for all $z \in X$, then $( \bar{x}, \bar{y}^{*}, \bar{\mu} ) \in \mathcal{S}(\hyperlink{CRD}{\mathrm{CRD}})$.
\end{itemize}
\end{thm}
\begin{proof}
According to Theorem \ref{Thm3-2}, there exist $y^{*} \in K^{+} \setminus \{ 0 \}$, $\mu := ( \mu_{1}, \mu_{2},\dots, \mu_{n} ) \in \mathbb{R}^{n}_{+}$, with $ \| y^{*} \| + \| \mu \| = 1$, and $\bar{v}_{i} \in \mathcal{V}_{i}$, $i = 1, 2,\dots, n$, satisfying
\begin{align}
  & 0 \in \medcup\limits_{ w^{*} \in \partial \langle y^{*}, f \rangle ( F ( \bar{x} ) ) } \partial \langle w^{*}, F \rangle ( \bar{x} ) + \sum_{i=1}^{n} \mu_{i} \, \mathrm{cl}^{*}\mathrm{co} \medPl \medcup \medBl \medcup\limits_{ v_{i}^{*} \in \partial_{x} g_{i} ( G_{i} ( \bar{x}, v_{i} ) ) } \partial_{x} \langle v_{i}^{*}, G_{i} \rangle ( \bar{x}, v_{i} ) \,\mid\, v_{i} \in \mathcal{V}_{i} ( \bar{x} ) \medBr \medPr, \nonumber \\
  \label{4-2}
  & \mu_{i} \, \max\limits_{ v_{i} \in \mathcal{V}_{i} } g_{i} ( G_{i} ( \bar{x}, v_{i} ) ) = 0, \quad i = 1, 2,\dots, n.
\end{align}
Letting $\bar{y}^{*} := y^{*}$ and $\bar{\mu} := ( \mu_{1}, \mu_{2},\dots, \mu_{n} )$, we get $( \bar{y}^{*}, \bar{\mu} ) \in K^{+} \setminus \{ 0 \} \times \mathbb{R}^{n}_{+}$. In addition, the inclusion $v_{i} \in \mathcal{V}_{i} ( \bar{x} )$ means that, for all $i \in \{ 1, 2,\dots, n \}$, one has $g_{i} ( G_{i} ( \bar{x}, v_{i} ) ) = \max\limits_{ u_{i} \in \mathcal{V}_{i} } g_{i} ( G_{i} ( \bar{x}, u_{i} ) )$. Hence, it follows from (\ref{4-2}) that $\mu_{i} \, g_{i} ( G_{i} ( \bar{x}, v_{i} ) ) = 0$, $i = 1, 2,\dots, n$. So $( \bar{x}, \bar{y}^{*}, \bar{\mu} ) \in C_{D}$.

(i) Let $( f \circ F, g \circ G )$ be type I pseudo convex at $z$ for all $z \in X$. For each $( z, y^{*}, \mu ) \in C_{D}$, employing Theorem \ref{Thm4-1}(i), we arrive at
\begin{equation*}
  \bar{f} ( \bar{x}, \bar{y}^{*}, \bar{\mu} ) = ( f \circ F ) ( \bar{x} ) \nprec \bar{f} ( z, y^{*}, \mu ).
\end{equation*}
Therefore $( \bar{x}, \bar{y}^{*}, \bar{\mu} ) \in \mathcal{S}^{w}(\hyperlink{CRD}{\mathrm{CRD}})$.

(ii) Let $( f \circ F, g \circ G )$ be type II pseudo convex at $z$ for all $z \in X$. Similarly, employing Theorem \ref{Thm4-1}(ii), we obtain
\begin{equation*}
  \bar{f} ( \bar{x}, \bar{y}^{*}, \bar{\mu} ) \npreceq \bar{f} ( z, y^{*}, \mu )
\end{equation*}
for each $( z, y^{*}, \mu ) \in C_{D}$. Thus $( \bar{x}, \bar{y}^{*}, \bar{\mu} ) \in \mathcal{S}(\hyperlink{CRD}{\mathrm{CRD}})$.
\end{proof}
\begin{rem}
\label{Rem4-1}
Theorem \ref{Thm4-1} and Theorem \ref{Thm4-2} improve both (\citealt[Theorem~4.1]{Saadati(2022)}) and (\citealt[Theorem~4.2]{Saadati(2022)}), both (\citealt[Theorem~5.2]{Fakhar(2018)}) and (\citealt[Corollary~5.4]{Fakhar(2018)}), and both (\citealt[Theorem~4.1]{Chuong(2016)}) and (\citealt[Theorem~4.3]{Chuong(2016)}).
\end{rem}
\begin{thm}
\label{Thm4-3}
\normalfont{\textbf{\textsc{(Strong Duality)}}}
Let $\bar{x} \in C$ be such that the robust $(\hyperlink{Def3-1}{\mathrm{KKT}})$ condition is satisfied at this point. Then there exists $( \bar{y}^{*}, \bar{\mu} ) \in K^{+} \setminus \{ 0 \} \times \mathbb{R}^{n}_{+}$ such that $( \bar{x}, \bar{y}^{*}, \bar{\mu} ) \in C_{D}$. Moreover,
\begin{itemize}
  \item[{\normalfont (i)}] If $( f \circ F, g \circ G )$ is type I pseudo convex at $z$ for all $z \in X$, then $( \bar{x}, \bar{y}^{*}, \bar{\mu} ) \in \mathcal{S}^{w}(\hyperlink{CRD}{\mathrm{CRD}})$ and $\bar{x} \in \mathcal{S}^{w}(\hyperlink{CRP}{\mathrm{CRP}})$.
  \item[{\normalfont (ii)}] If $( f \circ F, g \circ G )$ is type II pseudo convex at $z$ for all $z \in X$, then $( \bar{x}, \bar{y}^{*}, \bar{\mu} ) \in \mathcal{S}(\hyperlink{CRD}{\mathrm{CRD}})$ and $\bar{x} \in \mathcal{S}(\hyperlink{CRP}{\mathrm{CRP}})$.
\end{itemize}
\end{thm}
\begin{proof}
Using Definition \hyperlink{Def3-1}{3.1}, one can proceed similarly to the proof of Theorem \ref{Thm4-2}.
\end{proof}
\begin{rem}
\label{Rem4-2}
Theorem \ref{Thm4-1} and Theorem \ref{Thm4-3} develop both (\citealt[Theorem~4.1]{Saadati(2022)}) and (\citealt[Theorem~4.3]{Saadati(2022)}), and both (\citealt[Theorem~5.2]{Fakhar(2018)}) and (\citealt[Theorem~5.3]{Fakhar(2018)}).
\end{rem}

Finally, we establish converse duality relations between the problems $(\hyperlink{CRP}{\mathrm{CRP}})$ and $(\hyperlink{CRD}{\mathrm{CRD}})$.
\begin{thm}
\label{Thm4-4}
\normalfont{\textbf{\textsc{(Converse Duality)}}}
Let $( \bar{x}, \bar{y}^{*}, \bar{\mu} ) \in C_{D}$ be such that $\bar{x} \in C$.
\begin{itemize}
   \item[{\normalfont (i)}] If $( f \circ F, g \circ G )$ is type I pseudo convex at $\bar{x}$, then $\bar{x} \in \mathcal{S}^{w}(\hyperlink{CRP}{\mathrm{CRP}})$.
   \item[{\normalfont (ii)}] If $( f \circ F, g \circ G )$ is type II pseudo convex at $\bar{x}$, then $\bar{x} \in \mathcal{S}(\hyperlink{CRP}{\mathrm{CRP}})$.
\end{itemize}
\end{thm}
\begin{proof}
Since $( \bar{x}, \bar{y}^{*}, \bar{\mu} ) \in C_{D}$, we find $\bar{y}^{*} \in K^{+} \setminus \{ 0 \}$ and $\bar{\mu}_{i} \in \mathbb{R}^{n}_{+}$, $i = 1, 2,\dots, n$, such that
\begin{align}
  & 0 \in \medcup\limits_{ w^{*} \in \partial \langle \bar{y}^{*}, f \rangle ( F ( \bar{x} ) ) } \partial \langle w^{*}, F \rangle ( \bar{x} ) + \sum_{i=1}^{n} \bar{\mu}_{i} \, \mathrm{cl}^{*}\mathrm{co} \medPl \medcup \medBl \medcup\limits_{ v_{i}^{*} \in \partial_{x} g_{i} ( G_{i} ( \bar{x}, v_{i} ) ) } \partial_{x} \langle v_{i}^{*}, G_{i} \rangle ( \bar{x}, v_{i} ) \,\mid\, v_{i} \in \mathcal{V}_{i} ( \bar{x} ) \medBr \medPr, \nonumber \\
\label{4-3}
  & \bar{\mu}_{i} \, g_{i} ( G_{i} ( \bar{x}, v_{i} ) ) \ge 0, \,\, i = 1, 2,\dots, n.
\end{align}
On the other hand, since $\bar{x} \in C$, we have $( g_{i} \circ G_{i} ) ( \bar{x}, v_{i} ) \le 0$ for all $v_{i} \in \mathcal{V}_{i}$, $i = 1, 2,\dots, n$. Noting that $\bar{\mu}_{i} \, g_{i} ( G_{i} ( \bar{x}, v_{i} ) ) \le 0$ due to $\bar{\mu}_{i} \ge 0$. Hence, it follows from (\ref{4-3}) that $\bar{\mu}_{i} \, g_{i} ( G_{i} ( \bar{x}, v_{i} ) ) = 0$, and yields by taking into account Definition $\hyperlink{Def3-1}{3.1}$ that $\bar{x}$ is a $(\hyperlink{Def3-1}{\mathrm{KKT}})$ point of the problem $(\hyperlink{CRP}{\mathrm{CRP}})$. To finish the proof, it remains to apply Theorem \ref{Thm3-5}.
\end{proof}
\begin{rem}
\label{Rem4-3}
Theorem \ref{Thm4-4} reduces to (\citealt[Theorem~4.4]{Saadati(2022)}) for the case of problem (\hyperlink{UP}{UP}).
\end{rem}
\section{Applications}
\label{Sec5-Appl}
In this section, applications to necessary optimality conditions for an approximate uncertain multiobjective problem involving equality and inequality constraints and for a composite uncertain multiobjective problem with linear operators are given.

The class of multiobjective approximation problems plays a crucial role in optimization theory, as it provides a convenient framework that can effectively address a wide range of practical problems (see, e.g., \citealt{Gopfert(2003)B2, Jahn(2004)B}). Various problem domains, such as location problems, surrogate problems, inverse Stefan type problems, and optimal control problems, can be treated as  particular models of multiobjective approximation problems.

Let $X$ and $Y_{k}$, $k = 1, 2, \dots, l$, be Asplund spaces, let $r = ( r_{1}, r_{2}, \dots, r_{l} ): X \to \mathbb{R}^{l}$ be a vector-valued function, and let $T_{k} : X \to Y_{k}$, $k = 1, 2, \dots, l$, be linear operators. We select $y_{0}^{k} \in Y_{k}$, $\alpha_{k} \ge 0$, and $\beta_{k} \ge 1$, $k = 1, 2, \dots ,l$, and formulate the following \textit{approximate uncertain multiobjective optimization} problem
\begin{equation*}
\hypertarget{AUP}{}
  ( \mathrm{AUP} ) \qquad \min\nolimits_{\mathbb{R}^{l}_{+}} \,\,\, \medBl r (x) + ( \alpha_{1} \| T_{1} x - y_{0}^{1} \|^{\beta_{1}}, \alpha_{2} \| T_{2} x - y_{0}^{2} \|^{\beta_{2}}, \dots, \alpha_{l} \| T_{l} x - y_{0}^{l} \|^{\beta_{l}} ) \,\mid\, x \in C \medBr.
\end{equation*}
Here the feasible set $C$ is given by
\begin{align*}
  C := \medBl x \in X \,\mid\, &g_{i} ( x, v_{i} ) \le 0, \,\, i = 1, 2,\dots, n, \\
                               &h_{j} ( x, v_{n+j} ) = 0, \,\, j = 1, 2,\dots, m \medBr,
\end{align*}
where functions $g_{i}: X \times \mathcal{V}_{i} \to \mathbb{R}$ and $h_{j}: X \times \mathcal{V}_{n+j} \to \mathbb{R}$ define the constraints, $v_{i}$ and $v_{n+j}$ are uncertain parameters, and $v_{i} \in \mathcal{V}_{i}$ and $v_{n+j} \in \mathcal{V}_{n+j}$ for sequentially compact topological spaces. As mentioned before, we suppose that $r_{k}$, $k = 1, 2, \dots, l$, $g_{i}$, $i = 1, 2, \dots, n$, and $h_{j}$, $j = 1, 2, \dots, m$, are locally Lipschitzian at the point under consideration.

Now we are in a position to provide necessary conditions of the Fritz-John type for weakly robust efficient solutions of the problem $(\hyperlink{AUP}{\mathrm{AUP}})$, where the notions of such solutions are introduced similarly to the corresponding ones from Definition \hyperlink{Def1-1}{1.1}.
\begin{thm}
\label{Thm5-1}
Let $\bar{x}$ be a weakly robust efficient solution of the problem $(\hyperlink{AUP}{\mathrm{AUP}})$. Then there exist $\lambda_{k} \ge 0$, $k = 1, 2,\dots, l$, $\mu_{i} \ge 0$, $i = 1, 2,\dots, n$, $\sigma_{j} \ge 0$, $j = 1, 2,\dots, m$, not all zero, $y_{k}^{*} \in  Y_{k}^{*}$, $k = 1, 2,\dots, l$, with
\begin{equation*}
  \langle y_{k}^{*}, T_{k} \bar{x} - y_{0}^{k} \rangle = \| T_{k} \bar{x} - y_{0}^{k} \|^{\beta_{k}}, \quad
  \begin{cases}
    \| y^{*}_{k} \| \le 1 & \text{if } \beta_{k} = 1 \text{ and } T_{k} \bar{x} = y_{0}^{k}, \\
    \| y^{*}_{k} \| = \| T_{k} \bar{x} - y_{0}^{k} \|^{\beta_{k} - 1} & \text{otherwise,}
  \end{cases}
\end{equation*}
and $\bar{v}_{i} \in \mathcal{V}_{i}$, $i = 1, 2,\dots, n$, such that
\begin{align}
\label{5-1}
  & 0 \in \sum_{k=1}^{l} \lambda_{k} \partial r_{k} ( \bar{x} ) + \sum_{k=1}^{l} \lambda_{k} \alpha_{k} \beta_{k} T_{k}^{\top} y_{k}^{*} + \sum_{i=1}^{n} \mu_{i} \, \mathrm{cl}^{*}\mathrm{co} \medPl \medcup \medBl \partial_{x} g_{i} ( \bar{x}, v_{i} ) \,\mid\, v_{i} \in \mathcal{V}_{i} ( \bar{x} ) \medBr \medPr \nonumber \\
  & \hspace{2.8cm} + \sum_{j=1}^{m} \sigma_{j} \, \mathrm{cl}^{*}\mathrm{co} \medPl \medcup \medBl \partial_{x} h_{j} ( \bar{x}, v_{n+j} ) \medcup \partial_{x} ( - h_{j} ) ( \bar{x}, v_{n+j} ) \,\mid\, v_{n+j} \in \mathcal{V}_{n+j} \medBr \medPr, \\
  & \mu_{i} \, \max\limits_{ v_{i} \in \mathcal{V}_{i} } g_{i} ( \bar{x}, v_{i} ) = \mu_{i} \, g_{i} ( \bar{x}, \bar{v}_{i} ) = 0, \quad i = 1, 2,\dots, n. \nonumber
\end{align}
\end{thm}
\begin{proof}
For each $k = 1, 2, \dots, l$, put $\tilde{r}_{k} (x) := \alpha_{k} \| T_{k} x - y_{0}^{k} \|^{\beta_{k}}$, $x \in X$, and define vector-valued functions $F : X \to \mathbb{R}^{l} \times \mathbb{R}^{l}$, $f : \mathbb{R}^{l} \times \mathbb{R}^{l} \to \mathbb{R}^{l}$, $G_{i} : X \times \mathcal{V}_{i} \to X \times \mathcal{V}_{i}$, $i = 1, 2, \dots, n$, $G_{n+j} : X \times \mathcal{V}_{n+j} \to X \times \mathcal{V}_{n+j}$, $j = 1, 2, \dots, m$, $g_{i} : X \times \mathcal{V}_{i} \to \mathbb{R}$, $i = 1, 2, \dots, n$, and $g_{n+j} : X \times \mathcal{V}_{n+j} \to \mathbb{R}$, $j = 1, 2, \dots, m$, by
\begin{align*}
  F ( x ) & := ( r_{1} (x), r_{2} (x), \dots, r_{l} (x), \tilde{r}_{1} (x), \tilde{r}_{2} (x), \dots, \tilde{r}_{l} (x) ), \,\,\, x \in X, \\
  f ( w ) & := ( w_{1} + \tilde{w}_{1}, w_{2} + \tilde{w}_{2}, \dots, w_{l} + \tilde{w}_{l} ), \,\,\, w:= ( w_{1}, w_{2}, \dots, w_{l}, \tilde{w}_{1}, \tilde{w}_{2}, \dots, \tilde{w}_{l} ) \in \mathbb{R}^{l} \times \mathbb{R}^{l}, \\
  G_{i} ( x, v_{i} ) & := ( x, v_{i} ), \,\,\, x \in X, \,\, v_{i} \in \mathcal{V}_{i}, i = 1, 2, \dots, n, \\
  G_{n+j} ( x, v_{n+j} ) & := ( x, v_{n+j} ), \,\,\, x \in X, \,\, v_{n+j} \in \mathcal{V}_{n+j}, j = 1, 2, \dots, m, \\
  g ( z, u ) & := ( g_{1} ( z, u_{1} ), g_{2} ( z, u_{2} ), \dots, g_{n} ( z, u_{n} ), h_{1} ( z, u_{n+1} ), h_{2} ( z, u_{n+2} ), \dots, h_{m} ( z, u_{n+m} ) ), \\
  & \pushright{z \in X,\hspace{3.9cm}} \\
  & \pushright{u_{i} \in \mathcal{V}_{i}, \,\, i = 1, 2, \dots, n,\hspace{1.1cm}} \\
  & \pushright{u_{n+j} \in \mathcal{V}_{n+j}, \,\, j = 1, 2, \dots, m.}
\end{align*}
One can check that the problem $(\hyperlink{AUP}{\mathrm{AUP}})$ is a particular case of the problem $(\hyperlink{CUP}{\mathrm{CUP}})$ with $K := \mathbb{R}^{l}_{+}$. Employing Theorem \ref{Thm3-2} and proceeding as in the proof of Corollary \ref{Cor3-3}, we can find $y^{*} := ( \lambda_{1}, \lambda_{2}, \dots, \lambda_{l} ) \in K^{+} = \mathbb{R}^{l}_{+}$, $\mu := ( \mu_{1}, \mu_{2}, \dots, \mu_{n} ) \in \mathbb{R}^{n}_{+}$, $\sigma := ( \sigma_{1}, \sigma_{2},\dots, \sigma_{m} ) \in \mathbb{R}^{m}_{+}$, with $ \| y^{*} \| + \| \mu \| + \| \sigma \| \ne 0$, $\bar{v}_{i} \in \mathcal{V}_{i}$, $i = 1, 2, \dots, n$, and $\bar{v}_{n+j} \in \mathcal{V}_{n+j}$, $j = 1, 2, \dots, m$, satisfying
\begin{align}
  & 0 \in \medcup\limits_{ w^{*} \in \partial \langle y^{*}, f \rangle ( F ( \bar{x} ) ) } \partial \langle w^{*}, F \rangle ( \bar{x} ) + \sum_{i=1}^{n} \mu_{i} \, \mathrm{cl}^{*}\mathrm{co} \medPl \medcup \medBl \medcup\limits_{ v_{i}^{*} \in \partial_{x} g_{i} ( G_{i} ( \bar{x}, v_{i} ) ) } \partial_{x} \langle v_{i}^{*}, G_{i} \rangle ( \bar{x}, v_{i} ) \nonumber \\
  & \hspace{2mm} \,\mid\, v_{i} \in \mathcal{V}_{i} ( \bar{x} ) \medBr \medPr + \sum_{j=1}^{m} \sigma_{j} \, \mathrm{cl}^{*}\mathrm{co} \medPl \medcup \medBl \medcup\limits_{ v_{n+j}^{*} \in \partial_{x} h_{j} ( G_{n+j} ( \bar{x}, v_{n+j} ) ) \bigcup \partial_{x} - ( h_{j} ( G_{n+j} ( \bar{x}, v_{n+j} ) ) ) } \nonumber \\
\label{5-2}
  & \hspace{2mm} \partial_{x} \langle v_{n+j}^{*}, G_{n+j} \rangle ( \bar{x}, v_{n+j} ) \,\mid\, v_{n+j} \in \mathcal{V}_{n+j} ( \bar{x} ) \medBr \medPr \\
\intertext{and}
\label{5-3}
  &\left\{
    \begin{aligned}
      & \mu_{i} \, \max\limits_{ v_{i} \in \mathcal{V}_{i} } g_{i} ( G_{i} ( \bar{x}, v_{i} ) ) = \mu_{i} \, g_{i} ( G_{i} ( \bar{x}, \bar{v}_{i} ) ) = 0, \quad i = 1, 2, \dots, n, \\
      & \sigma_{j} \, \max\limits_{ v_{n+j} \in \mathcal{V}_{n+j} } h_{j} ( G_{n+j} ( \bar{x}, v_{n+j} ) ) = \sigma_{j} \, h_{j} ( G_{n+j} ( \bar{x}, \bar{v}_{n+j} ) ) = 0, \quad j = 1, 2, \dots, m.
    \end{aligned}
  \right.
\end{align}
It follows from the above definition that for each $w \in \mathbb{R}^{l} \times \mathbb{R}^{l}$ one has
\begin{equation*}
  \langle y^{*}, f \rangle (w) = \sum_{k=1}^{l} \lambda_{k} ( w_{k} + \tilde{w}_{k} ),
\end{equation*}
and hence
\begin{equation*}
  \partial \langle y^{*}, f \rangle (w) = \{ ( \lambda_{1}, \lambda_{2}, \dots, \lambda_{l}, \lambda_{1}, \lambda_{2}, \dots, \lambda_{l} ) \}.
\end{equation*}
This shows that
\begin{equation*}
  \partial \langle y^{*}, f \rangle ( F ( \bar{x} ) ) = \{ w^{*} := ( \lambda_{1}, \lambda_{2}, \dots, \lambda_{l}, \lambda_{1}, \lambda_{2}, \dots, \lambda_{l} ) \}.
\end{equation*}
On the other side, it holds by the definition that
\begin{equation*}
  \langle w^{*}, F \rangle ( x ) = \sum_{k=1}^{l} \lambda_{k} r_{k} ( x ) + \sum_{k=1}^{l} \lambda_{k} \tilde{r}_{k} ( x ), \quad x \in X.
\end{equation*}
Thus applying the sum rule of Lemma \ref{Lem2-3}, we get
\begin{equation*}
  \partial \langle w^{*}, F \rangle ( \bar{x} ) = \partial \Big( \sum_{k=1}^{l} \lambda_{k} r_{k} + \sum_{k=1}^{l} \lambda_{k} \tilde{r}_{k} \Big) ( \bar{x} ) \subset \sum_{k=1}^{l} \lambda_{k} \partial r_{k} ( \bar{x} ) + \sum_{k=1}^{l} \lambda_{k} \partial \tilde{r}_{k} ( \bar{x} ).
\end{equation*}
In addition putting
\begin{equation*}
  p_{k} ( x ) := T_{k} x - y^{k}, \quad x \in X, \qquad q_{k} ( y ) := \alpha_{k} \| y \|^{\beta_{k}}, \quad y \in Y_{k},
\end{equation*}
we have from Lemmas \ref{Lem2-2} and \ref{Lem2-6} that
\begin{equation*}
  \partial \tilde{r}_{k} ( \bar{x} ) = \partial ( p_{k} \circ q_{k} ) ( \bar{x} ) \subset \big\{ \alpha_{k} \beta_{k} T_{k}^{\top} y_{k}^{*} \,\mid\, y_{k}^{*} \in Y_{k}^{*} \big\}
\end{equation*}
with
\begin{equation}
\label{5-4}
  \begin{cases}
    \| y^{*}_{k} \| \le 1 & \text{if } \beta_{k} = 1 \text{ and } T_{k} \bar{x} = y_{0}^{k}, \\
    \langle y_{k}^{*}, T_{k} \bar{x} - y_{0}^{k} \rangle = \| T_{k} \bar{x} - y_{0}^{k} \|^{\beta_{k}} \text{ and } \| y^{*}_{k} \| = \| T_{k} \bar{x} - y_{0}^{k} \|^{\beta_{k} - 1} & \text{otherwise.}
  \end{cases}
\end{equation}
So, we arrive at the following inclusion
\begin{equation}
\label{5-5}
  \medcup\limits_{ w^{*} \in \partial \langle y^{*}, f \rangle ( F ( \bar{x} ) ) } \partial \langle w^{*}, F \rangle ( \bar{x} ) \subset \sum_{k=1}^{l} \lambda_{k} \partial r_{k} ( \bar{x} ) + \Bigg\{ \sum_{k=1}^{l} \lambda_{k} \alpha_{k} \beta_{k} T_{k}^{\top} y_{k}^{*} \,\mid\, y_{k}^{*} \in Y_{k}^{*} \Bigg\},
\end{equation}
where $y_{k}^{*}$ for $k = 1, 2,\dots,l$ satisfy (\ref{5-4}). Similarly to the above, we obtain
\begin{align}
\label{5-6}
  & \medcup\limits_{ v_{i}^{*} \in \partial_{x} g_{i} ( G_{i} ( \bar{x}, v_{i} ) ) } \partial_{x} \langle v_{i}^{*}, G_{i} \rangle ( \bar{x}, v_{i} ) = \partial_{x} g_{i} ( \bar{x}, v_{i} ), \quad i = 1, 2,\dots, n,  \\
   \intertext{and}
\label{5-7}
  & \medcup\limits_{ v_{n+j}^{*} \in \partial_{x} h_{j} ( G_{n+j} ( \bar{x}, v_{n+j} ) ) \bigcup \partial_{x} - ( h_{j} ( G_{n+j} ( \bar{x}, v_{n+j} ) ) ) } \partial_{x} \langle v_{n+j}^{*}, G_{n+j} \rangle ( \bar{x}, v_{n+j} ) \nonumber \\
  \subset \,\, & \partial_{x} h_{j} ( \bar{x}, v_{n+j} ) \medcup \partial_{x} ( - h_{j} ) ( \bar{x}, v_{n+j} ), \quad j = 1, 2,\dots, m.
\end{align}
Combining (\ref{5-2}) with (\ref{5-5})-(\ref{5-7}), we finally get (\ref{5-1}). Note that in our setting (\ref{5-3}) collapses to the following
\begin{equation*}
  \mu_{i} \, \max\limits_{ v_{i} \in \mathcal{V}_{i} } g_{i} ( \bar{x}, v_{i} ) = \mu_{i} \, g_{i} ( \bar{x}, \bar{v}_{i} ) = 0
\end{equation*}
for each $i = 1, 2,\dots, n$. This justifies the last statement of the theorem and completes the proof.
\end{proof}
\begin{exa}
\label{Exa5-1}
Let $X := \mathbb{R}^{2}$, $Y_{k} := \mathbb{R}^{2}$, $k = 1, 2, 3$, $\mathcal{V}_{i} := [-1, -\dfrac{1}{4}]$, $i = 1, 2$, and $\mathcal{V}_{2+j} := [-1, -\dfrac{1}{4}]$, $j = 1, 2$. Take the following approximate uncertain multiobjective optimization problem
\begin{equation*}
\hypertarget{AUPexa5.1}{}
  \begin{aligned}
    ( \mathrm{AUP} ) \qquad \min\nolimits_{\mathbb{R}^{3}_{+}} \,\,\, \medBl &r (x) + ( \alpha_{1} \| T_{1} x - y_{0}^{1} \|^{\beta_{1}}, \alpha_{2} \| T_{2} x - y_{0}^{2} \|^{\beta_{2}}, \alpha_{3} \| T_{3} x - y_{0}^{3} \|^{\beta_{3}} ) \,\mid\, x := ( x_{1}, x_{2} ) \in X, \\
    &g_{i} ( x, v_{i} ) \le 0, \,\, i = 1, 2 \quad \text{and} \quad h_{j} ( x, v_{2+j} ) = 0, \,\, j = 1,2 \medBr,
  \end{aligned}
\end{equation*}
where $r : X \to \mathbb{R}^{3}$, $r := (r_{1}, r_{2}, r_{3})$, is defined by
\begin{equation*}
  \left\{\begin{aligned}
           r_{1}(x_{1},x_{2}) &:= 3 |x_{1}| + \frac{2}{5} x_{2} + \frac{4}{5}, \\
           r_{2}(x_{1},x_{2}) &:= \dfrac{1}{4} x_{1}^{2} + 2, \\
           r_{3}(x_{1},x_{2}) &:= 2 |x_{1}| - \frac{1}{8} x_{2}^{2} + 1,
         \end{aligned}
  \right.
\end{equation*}
$g_{i} : X \times \mathcal{V}_{i} \to \mathbb{R}$, $i = 1, 2$, and $h_{j} : X \times \mathcal{V}_{2+j} \to \mathbb{R}$, $j = 1, 2$, are given respectively by
\begin{equation*}
  \left\{\begin{aligned}
           g_{1}(x_{1},x_{2}, v_{1}) &:= \frac{1}{4} v_{1}^{2} |x_{1}| + \frac{1}{2} v_{1}^{2} x_{2} + \frac{1}{4} |v_{1}|, \\
           g_{2}(x_{1},x_{2}, v_{2}) &:= \frac{1}{8} x_{1}^{2} + |v_{2}| x_{2} + |v_{2}| + \frac{1}{4},
         \end{aligned}
  \right.
  \hspace{5mm}
  \left\{\begin{aligned}
           h_{1}(x_{1},x_{2}, v_{3}) &:= v_{3} ( - 3 x_{1} + x_{2} + 2 ), \\
           h_{2}(x_{1},x_{2}, v_{4}) &:= v_{4} ( - 3 x_{1} - x_{2} - 2 ),
         \end{aligned}
  \right.
\end{equation*}
$T_{1} := \begin{pmatrix} 0 & \dfrac{1}{2} \\ 1 & 0 \end{pmatrix}$, $T_{2} := \begin{pmatrix} 1 & 0 \\ 0 & \dfrac{1}{2} \end{pmatrix}$, $T_{3} := \begin{pmatrix} 0 & \dfrac{1}{2} \\ 0 & \dfrac{1}{2} \end{pmatrix}$, $y_{0}^{1} := ( -1, 0 ) \in Y_{1}$, $y_{0}^{2} := (0, -1)\in Y_{2}$, $y_{0}^{3} := ( -1, -1 ) \in Y_{3}$, $\alpha_{1} := 2$, $\alpha_{2} := \alpha_{3} := 1$, and $\beta_{1} := \beta_{2} := \beta_{3} := 1$. Note that from the constraints of inequality and equality one get
\begin{equation*}
  C = \medBl ( x_{1}, x_{2} ) \in X \,\mid\, - 20 \le x_{2} \le - 2, \quad \text{and} \quad x_{2} = 3 x_{1} - 2, \,\, x_{2} = - 3 x_{1} - 2 \medBr.
\end{equation*}

Now, suppose that $\bar{x} := ( 0, -2 ) \in C$ and $x := ( x_{1}, x_{2} ) \in C$. It is easy to check that $\bar{x}$ is a weakly robust efficient solution of the problem $(\hyperlink{AUPexa5.1}{\mathrm{AUP}})$. One can also see that
\begin{align*}
  \phi_{1}(\bar{x}) &= \max\limits_{v_{1} \in \mathcal{V}_{1}} g_{1}(\bar{x}, v_{1}) = \max\limits_{v_{1} \in \mathcal{V}_{1}} ( -v_{1}^{2} + \dfrac{1}{4} | v_{1} | ) = 0, \\
  \phi_{2}(\bar{x}) &= \max\limits_{v_{2} \in \mathcal{V}_{2}} g_{2}(\bar{x}, v_{2}) = \max\limits_{v_{2} \in \mathcal{V}_{2}} (- |v_{2}| + \dfrac{1}{4}) = 0.
\end{align*}
Hence $\phi(\bar{x}) = \max \big\{ \phi_{1}(\bar{x}), \phi_{2}(\bar{x}) \big\} = 0$, $\mathcal{V}_{1}(\bar{x}) = \{-\dfrac{1}{4}\}$, and $\mathcal{V}_{2}(\bar{x}) = \{-\dfrac{1}{4}\}$. Performing elementary calculations gives us from the definition that $\partial r_{1}(\bar{x}) = [-3, 3] \times \{\dfrac{2}{5}\}$, $\partial r_{2}(\bar{x}) = ( 0, 0 )$, $\partial r_{3}(\bar{x}) = [-2, 2] \times \{\dfrac{1}{2}\}$,
\begin{align*}
  &\mathrm{cl}^{*}\mathrm{co} \medPl \partial_{x}g_{1} ( \bar{x}, v_{1} = -\dfrac{1}{4} ) \medPr = [ -\dfrac{1}{64}, \dfrac{1}{64} ] \times \{\dfrac{1}{32}\}, \quad \mathrm{cl}^{*}\mathrm{co} \medPl \partial_{x}g_{2} ( \bar{x}, v_{2} = -\dfrac{1}{4} ) \medPr = (0, \dfrac{1}{4}),\\
  \intertext{and}
  &\mathrm{cl}^{*}\mathrm{co}  \medPl \medcup \medBl \partial_{x} h_{1} ( \bar{x}, v_{3} ) \medcup \partial_{x} ( - h_{1} ) ( \bar{x}, v_{3} ) \,\mid\, v_{3} \in \mathcal{V}_{3} \medBr \medPr = \big\{ ( v_{3}, - \dfrac{1}{3} v_{3} ) \,\mid\, v_{3} \in [ -3, 3 ]  \big\},\\
  &\mathrm{cl}^{*}\mathrm{co}  \medPl \medcup \medBl \partial_{x} h_{2} ( \bar{x}, v_{4} ) \medcup \partial_{x} ( - h_{2} ) ( \bar{x}, v_{4} ) \,\mid\, v_{4} \in \mathcal{V}_{4} \medBr \medPr = \big\{ ( v_{4}, \dfrac{1}{3} v_{4} ) \,\mid\, v_{4} \in [ -3, 3 ]  \big\}.
\end{align*}
Finally, observe that there exist $\lambda = ( 1, 0, 1) \in \mathbb{R}^{3}_{+}$, $\mu = ( 0, 1 ) \in \mathbb{R}^{2}_{+}$, $\sigma = ( 1 , 0 ) \in \mathbb{R}^{2}_{+}$, $y^{*}_{1} = ( -\dfrac{2}{5}, 0 ) \in Y^{*}_{1}$, $y^{*}_{2} = ( 0, 0 ) \in Y^{*}_{2}$, and $y^{*}_{3} = ( -\dfrac{1}{2}, -\dfrac{1}{2} ) \in Y^{*}_{3}$ satisfying
\begin{align*}
  0 &=
  \begin{pmatrix}
    1 & 0  & 1
  \end{pmatrix}
  \begin{pmatrix}
    -\dfrac{3}{4} & 0 & 0 \\
    \dfrac{2}{5} & 0 & \dfrac{1}{2}
  \end{pmatrix}
  +
  \begin{pmatrix}
    0 & 1
  \end{pmatrix}
  \begin{pmatrix}
    0 & 0 \\
    \dfrac{1}{32} & \dfrac{1}{4}
  \end{pmatrix}
  +
  \begin{pmatrix}
    1 & 0
  \end{pmatrix}
  \begin{pmatrix}
    \dfrac{3}{4} & 0 \\
    -\dfrac{1}{4} & 0
  \end{pmatrix}
  \\
  &+ 2
  \begin{pmatrix}
  0 & 1 \\
  \dfrac{1}{2} & 0 \\
  \end{pmatrix}
  \begin{pmatrix}
  - \dfrac{2}{5} \\
  0
  \end{pmatrix}
  + 0
  \begin{pmatrix}
  1 & 0 \\
  0 & \dfrac{1}{2} \\
  \end{pmatrix}
  \begin{pmatrix}
  0 \\
  0
  \end{pmatrix}
  +
  \begin{pmatrix}
  0 & 0 \\
  \dfrac{1}{2} & \dfrac{1}{2}
  \end{pmatrix}
  \begin{pmatrix}
  - \dfrac{1}{2} \\
  - \dfrac{1}{2}
  \end{pmatrix},
\end{align*}
and $\mu_{i} \, \max\limits_{v_{i} \in \mathcal{V}_{i}} g_{i}(\bar{x}, v_{i}) = 0$ for $i = 1, 2$.
\end{exa}

Next let $T_{k} : X \to Y_{k}$, $k = 1, 2,\dots, l$ be linear operators and $f_{k} : Y_{k} \to \mathbb{R}$, $k = 1, 2,\dots, l$ be local Lipschitz functions between Asplund spaces. We consider the following composite uncertain multiobjective optimization problem with \textit{linear} operators
\begin{equation*}
\hypertarget{CUL}{}
  ( \mathrm{CUL} ) \qquad \min\nolimits_{\mathbb{R}^{l}_{+}} \,\,\, \medBl ( f_{1} ( T_{1} x ), f_{2} ( T_{2} x ), \dots, f_{l} ( T_{l} x ) ) \,\mid\, x \in C \medBr.
\end{equation*}
Here the feasible set $C$ is given by
\begin{align*}
  C := \medBl x \in X \,\mid\, g_{i} ( x, v_{i} ) \le 0, \,\, i = 1, 2,\dots, n \medBr,
\end{align*}
where function $g_{i}: X \times \mathcal{V}_{i} \to \mathbb{R}$ is locally Lipschitz, $v_{i}$ is uncertain parameter, and $v_{i} \in \mathcal{V}_{i}$ for sequentially compact topological spaces.

We now derive necessary conditions of the Fritz-John form for weakly robust efficient solutions of the problem $(\hyperlink{CUL}{\mathrm{CUL}})$.
\begin{thm}
\label{Thm5-2}
Let $\bar{x}$ be a weakly robust efficient solution of the problem $(\hyperlink{CUL}{\mathrm{CUL}})$. Then there exist $\lambda_{k} \ge 0$, $k = 1, 2,\dots, l$, $\mu_{i} \ge 0$, $i = 1, 2,\dots, n$, with $\| ( \lambda_{1}, \lambda_{2},\dots, \lambda_{l} ) \| + \| ( \mu_{1}, \mu_{2},\dots, \mu_{n} ) \| = 1$, and $\bar{v}_{i} \in \mathcal{V}_{i}$, $i = 1, 2,\dots, n$, such that
\begin{align}
\label{5-8}
  & 0 \in \sum\limits_{k=1}^{l} \lambda_{k} T_{k}^{\top} \partial f_{k} ( \bar{y}_{k} ) + \sum_{i=1}^{n} \mu_{i} \, \mathrm{cl}^{*}\mathrm{co} \medPl \medcup \medBl \partial_{x} g_{i} ( \bar{x}, v_{i} ) \,\mid\, v_{i} \in \mathcal{V}_{i} ( \bar{x} ) \medBr \medPr, \\
\label{5-9}
  & \mu_{i} \, \max\limits_{ v_{i} \in \mathcal{V}_{i} } g_{i} ( \bar{x}, v_{i} ) = \mu_{i} \, g_{i} ( \bar{x}, \bar{v}_{i} ) = 0, \quad i = 1, 2,\dots, n,
\end{align}
where $\bar{y}_{k} = T_{k} \bar{x}$ for $k = 1, 2,\dots, l$.
\end{thm}
\begin{proof}
We first observe that the problem $(\hyperlink{CUL}{\mathrm{CUL}})$ is a special case of the composite uncertain multiobjective optimization problem $(\hyperlink{CUP}{\mathrm{CUP}})$, where $Z = X$, $G = ( G_{1}, G_{2},\dots, G_{n} )$ is an identical map, the functions $F : X \to W := Y_{1} \times Y_{2} \times \dots \times Y_{l}$ and $f : W \to Y := \mathbb{R}^{l}$ are defined, respectively, by
\begin{align*}
  F (x) & := ( T_{1} x, T_{2} x, \dots, T_{l} x ), \quad x \in X, \\
  f (w) & := ( f_{1} ( w_{1} ), f_{2} ( w_{2} ),\dots, f_{l} ( w_{l} ) ), \quad w := ( w_{1}, w_{2},\dots, w_{l} ) \in Y_{1} \times Y_{2} \times \dots \times Y_{l},
\end{align*}
and $K := \mathbb{R}^{l}_{+}$. Applying Theorem \ref{Thm3-2}, we find $y^{*} := ( \lambda_{1}, \lambda_{2}, \dots, \lambda_{l} ) \in K^{+} = \mathbb{R}^{l}_{+}$, $\mu := ( \mu_{1}, \mu_{2}, \dots, \mu_{n} ) \in \mathbb{R}^{n}_{+}$, with $ \| y^{*} \| + \| \mu \| = 1$, and $\bar{v}_{i} \in \mathcal{V}_{i}$, $i = 1, 2, \dots, n$, such that
\begin{align}
\label{5-10}
  & 0 \in \medcup\limits_{ w^{*} \in \partial \langle y^{*}, f \rangle ( F ( \bar{x} ) ) } \partial \langle w^{*}, F \rangle ( \bar{x} ) + \sum_{i=1}^{n} \mu_{i} \, \mathrm{cl}^{*}\mathrm{co} \medPl \medcup \medBl \medcup\limits_{ v_{i}^{*} \in \partial_{x} g_{i} ( G_{i} ( \bar{x}, v_{i} ) ) } \partial_{x} \langle v_{i}^{*}, G_{i} \rangle ( \bar{x}, v_{i} ) \,\mid\, v_{i} \in \mathcal{V}_{i} ( \bar{x} ) \medBr \medPr, \\
\label{5-11}
  & \mu_{i} \, \max\limits_{ v_{i} \in \mathcal{V}_{i} } g_{i} ( G_{i} ( \bar{x}, v_{i} ) ) = \mu_{i} \, g_{i} ( G_{i} ( \bar{x}, \bar{v}_{i} ) ) = 0, \quad i = 1, 2,\dots, n.
\end{align}
For each $k = 1, 2,\dots, l$, consider a mapping $\Psi_{k} : Y_{1} \times Y_{2} \times \dots \times Y_{l} \to Y_{k}$ given by $\Psi_{k} ( w ) := w_{k}$, $w := ( w_{1}, w_{2},\dots, w_{l} ) \in Y_{1} \times Y_{2} \times \dots \times Y_{l}$. Using the definition, we have
\begin{equation*}
\langle y^{*}, f \rangle ( w ) = \sum_{k=1}^{l} \lambda_{k} ( f_{k} \circ \Psi_{k}) (w), \quad w \in Y_{1} \times Y_{2} \times \dots \times Y_{l},
\end{equation*}
and hence invoking Lemmas \ref{Lem2-2} and \ref{Lem2-3},
\begin{align*}
\partial \langle y^{*}, f \rangle ( w ) &\subset \sum_{k=1}^{l} \lambda_{k} ( f_{k} \circ \Psi_{k}) (w) \\
&\subset \lambda_{1} \partial f_{1} ( w_{1} ) \times \lambda_{2} \partial f_{2} ( w_{2} ) \times \dots \times \lambda_{l} \partial f_{l} ( w_{l} ).
\end{align*}
This gives
\begin{equation}
\label{5-12}
\partial \langle y^{*}, f \rangle ( F ( \bar{x} ) ) \subset \lambda_{1} \partial f_{1} ( y_{1} ) \times \lambda_{2} \partial f_{2} ( y_{2} ) \times \dots \times \lambda_{l} \partial f_{l} ( y_{l} ),
\end{equation}
where $y_{k} = T_{k} x$ for $k = 1, 2,\dots, l$.

Now picking any $w^{*} \in \partial \langle y^{*}, f \rangle ( F ( \bar{x} ) )$ and taking (\ref{5-12}) into account, there exist $w^{*}_{k} \in \partial f_{k} ( \bar{y}_{k} )$, $k = 1, 2,\dots, l$, with $\bar{y}_{k} = T_{k} \bar{x}$ satisfying $w^{*} = ( \lambda_{1} w_{1}^{*^{\top}}, \lambda_{2} w_{2}^{*^{\top}},\dots, \lambda_{l} w_{l}^{*^{\top}} )$. Then it follows from the latter that
\begin{equation*}
\langle w^{*}, F \rangle (x) = \sum_{k=1}^{l} \lambda_{k} w_{k}^{*^{\top}} T_{k} x, \quad x \in X,
\end{equation*}
and thus $\partial \langle w^{*}, F \rangle ( \bar{x} ) = \sum\limits_{k=1}^{l} \lambda_{k} T_{k}^{\top} w_{k}^{*}$. So, we get the inclusion
\begin{equation}
\label{5-13}
\medcup\limits_{ w^{*} \in \partial \langle y^{*}, f \rangle ( F ( \bar{x} ) ) } \partial \langle w^{*}, F \rangle ( \bar{x} ) \subset \sum\limits_{k=1}^{l} \lambda_{k} T_{k}^{\top} \partial f_{k} ( \bar{y}_{k} ).
\end{equation}
Similarly, we arrive at
\begin{equation}
\label{5-14}
\medcup\limits_{ v_{i}^{*} \in \partial_{x} g_{i} ( G_{i} ( \bar{x}, v_{i} ) ) } \partial_{x} \langle v_{i}^{*}, G_{i} \rangle ( \bar{x}, v_{i} ) = \partial_{x} g_{i} ( \bar{x}, v_{i} ), \quad i = 1, 2,\dots, n,
\end{equation}
since $G$ is an identical map. Finally, combining (\ref{5-10}) with the relations in (\ref{5-13}) and (\ref{5-14}), we obtain (\ref{5-8}). Furthermore, in our setting (\ref{5-11}) reduces to (\ref{5-9}), and so the proof of the theorem is complete.
\end{proof}
\begin{exa}
\label{Exa5-2}
Let $X$, $Y_{k}$, $k = 1, 2, 3$, $\mathcal{V}_{i}$, $i = 1, 2$, $T_{k}$, $k = 1, 2, 3$, $r_{k}$, $k = 1, 2, 3$, and $g_{i}$, $i = 1, 2$, be the same as Example \ref{Exa5-1}, and let $f_{k} := r_{k}$, $k = 1, 2, 3$. Take the following composite uncertain multiobjective optimization problem with linear operators
\begin{equation*}
\hypertarget{CULexa5.2}{}
  \begin{aligned}
    ( \mathrm{CUL} ) \qquad \min\nolimits_{\mathbb{R}^{3}_{+}} \,\,\, \medBl &( f_{1} ( T_{1}x ), f_{2} ( T_{2}x ), f_{3} ( T_{3}x ) ) \,\mid\, x := ( x_{1}, x_{2} ) \in X, \,\, g_{i} ( x, v_{i} ) \le 0, \,\, i = 1, 2 \medBr.
  \end{aligned}
\end{equation*}
Note that the robust feasible set is given by
\begin{equation*}
  C = \medBl ( x_{1}, x_{2} ) \in X \,\mid\, |x_{1}| \le 1 \text{ and } x_{2} \le - \dfrac{1}{2} |x_{1}| - 2 \medBr \medcup \medBl ( x_{1}, x_{2} ) \in X \,\mid\, |x_{1}| > 1 \text{ and } x_{2} \le - \dfrac{1}{2} x_{1}^{2} - 2 \medBr.
\end{equation*}
Suppose that $ \bar{x} := ( 0, -2 ) \in C$. It is easy to check that $f_{2} ( T_{2} x ) - f_{2} ( T_{2} \bar{x} ) \ge 0$ for all $x \in C$, i.e., $\bar{x}$ is a weakly robust efficient solution of the problem $(\hyperlink{CULexa5.2}{\mathrm{CUL}})$. Taking into account that $\bar{y}_{1} = T_{1} \bar{x} = \begin{pmatrix} -1 \\ 0 \end{pmatrix}$, $\bar{y}_{2} = T_{2} \bar{x} = \begin{pmatrix} 0 \\ -1 \end{pmatrix}$, and $\bar{y}_{3} = T_{3} \bar{x} = \begin{pmatrix} -1 \\ -1 \end{pmatrix}$, we have $\partial r_{1}(\bar{y}) = [-3, 3] \times \{\dfrac{2}{5}\}$, $\partial r_{2}(\bar{y}) = ( 0, 0 )$, and $\partial r_{3}(\bar{y}) = [-2, 2] \times \{\dfrac{1}{4}\}$.
Finally, we can find $\lambda = ( 0, 0, \dfrac{1}{2} ) \in \mathbb{R}^{3}_{+}$ and $\mu = ( 0, \dfrac{1}{2} ) \in \mathbb{R}^{2}_{+}$, with $\| \lambda \| + \| \mu \| = 1$, such that
\begin{align*}
  0 &= 0
  \begin{pmatrix}
  0 & 1 \\
  \dfrac{1}{2} & 0 \\
  \end{pmatrix}
  \begin{pmatrix}
  0 \\
  \dfrac{2}{5}
  \end{pmatrix}
  + 0
  \begin{pmatrix}
  1 & 0 \\
  0 & \dfrac{1}{2}
  \end{pmatrix}
  \begin{pmatrix}
  0 \\
  0
  \end{pmatrix}
  + \dfrac{1}{2}
  \begin{pmatrix}
  0 & 0 \\
  \dfrac{1}{2} & \dfrac{1}{2}
  \end{pmatrix}
  \begin{pmatrix}
  -\dfrac{3}{4} \\
  \dfrac{1}{4}
  \end{pmatrix}
  +
  \begin{pmatrix}
    0 & \dfrac{1}{2}
  \end{pmatrix}
  \begin{pmatrix}
    0 & 0 \\
    \dfrac{1}{32} & \dfrac{1}{4}
  \end{pmatrix}.
\end{align*}
\end{exa}
\pdfbookmark[section]{References}{bibliography}
\bibliographystyle{apa}
\small
\bibliography{mybib.bib}
\end{document}